\documentclass[11pt]{article}
\usepackage{amssymb,amsthm,amsmath,amsfonts}
\usepackage[latin1]{inputenc}
\usepackage{graphicx}
\usepackage{verbatim}
\usepackage{array}
\usepackage[colorlinks=true,urlcolor=blue,
citecolor=red,linkcolor=blue,linktocpage,pdfpagelabels,
bookmarksnumbered,bookmarksopen]{hyperref}
\usepackage{color}
\newtheorem{Theorem}{Theorem}[section]

\newtheorem{theorem}[Theorem]{Theorem}

\newtheorem{remark}[Theorem]{Remark}

\newtheorem{proposition}[Theorem]{Proposition}

\newtheorem{lemma}[Theorem]{Lemma}

\newtheorem{definition}[Theorem]{Definition}

\newtheorem{corollary}[Theorem]{Corollary}

\definecolor{red}{rgb}{0.9,0,0}

\definecolor{azul}{rgb}{0,0,1}

\def\R {{\rm I}\hskip -0.85mm{\rm R}}
\renewcommand{\O}{\Omega}
\def\D{\Delta}
\def\L{\Lambda}
\renewcommand{\l}{\lambda}
\def\m{\mu}

\def\s{\sigma}
\def\v{\varphi}

\def\g{\gamma}
\def\n{\nabla}

\def\bu{{\bf u}}

\title{Generalized eigenvalue problem for an interface elliptic equation}

\author{Braulio B. V. Maia$^1$, M\'onica Molina-Becerra$^2$, Cristian Morales-Rodrigo$^3$, Antonio Su\'arez$^3$}
\date{}

\begin{document}
	
	\pretolerance10000
	
	\maketitle
	
	1. Universidade Federal Rural da Amazonia, Campus de Capitao-Poco, PA, Brazil.
	
	2. Dpto. Matem\'atica Aplicada II, Escuela Polit\'ecnica Superior, C. Virgen de Africa, 7, 41011, Univ. de Sevilla, Sevilla, Spain.
	
	3. Dpto Ecuaciones Diferenciales y An\'alisis Num\'erico and IMUS, Fac. de Matem\'aticas, Univ. de Sevilla, Sevilla, Spain.
	\vskip 0.5cm
	
	e-mails: braulio.maia@ufra.edu.br, monica@us.es, cristianm@us.es, suarez@us.es.

	\begin{abstract}
		In this paper we deal with an eigenvalue problem in an interface elliptic equation. We characterize the set of principal eigenvalues as a level set of a concave and regular function. As application, we study a problem arising in population dynamics. In these problems each species lives in a subdomain, and they interact in a common border, which acts as a geographical barrier. 
	\end{abstract}

	\noindent{\bf\small Keywords:} {\small interface, principal eigenvalue.}
	
	\noindent{\bf\small MSC2010:} {\small 35A15, 35B33, 35B25, and 35J60.}
	
	\section{Introduction}

\begin{figure}
\centering
\includegraphics[scale=0.6]{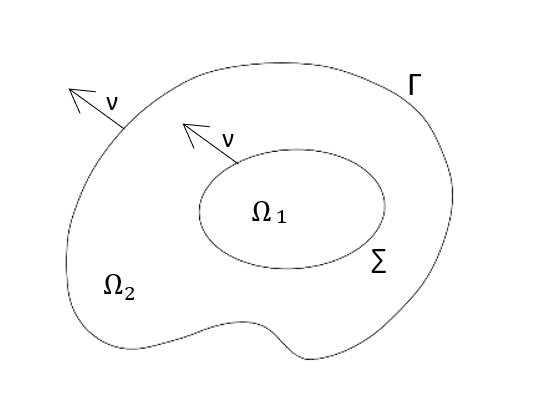}
\caption{A possible configuration of the domain $\O=\O_1\cup \O_2 \cup \Sigma$.}
\label{figure1}
\end{figure}

Recently, the following  semilinear interface problems have been analyzed
\begin{equation}
\label{sisintro}
 \left\{\begin{array}{ll}
-\D u_i=\l f_i(x,u_i) & \text{in $\Omega_i$, $i=1,2$},\\
\noalign{\smallskip}
\partial_{\nu}u_i=\gamma_i(u_2-u_1)& \text{on $\Sigma$},\\
\partial_{\nu}u_2=0& \text{on $\Gamma$,}
\end{array}
\right. 
\end{equation}
where $\O$ is a bounded domain of $\R^N$ with
$$
\O=\O_1\cup \O_2 \cup \Sigma,
$$
with $\O_i$ subdomains, with internal interface $\Sigma=\partial \O_1$, and $\Gamma=\partial\O_2\setminus\Sigma$, $\nu_i$ is the outward normal to $\O_i$, and we call $\nu:=\nu_1=-\nu_2$ (see Figure \ref{figure1} where we have illustrated an example of $\O$).   

In (\ref{sisintro}), $u_i$ represents the density of a species inhabiting in $\O_i$, and they interact on $\Sigma$ under the so called Kedem-Katchalsky conditions (see \cite{KedemK}), and it means that the flux is proportional to the jump of the function through $\Sigma$ (see  \cite{chenCPDE}, \cite{chen2}, \cite{cp}, \cite{cia2} and references therein). Here, $f_i:\O_i\times \R\mapsto \R$ are regular functions, $\g_i>0$ stands for the proportional coefficient of the jump and $1/\l$, $\l>0$, is a real parameter  representing the diffusion coefficient of the species, the same in both subdomains. It seems natural to consider two different diffusion coefficients, one for each species, that is, a problem as  
\begin{equation}
\label{sisintro2}
 \left\{\begin{array}{ll}
-\D u_i=\l_i f_i(x,u_i) & \text{in $\Omega_i$},\\
\noalign{\smallskip}
\partial_{\nu}u_i=\gamma_i(u_2-u_1)& \text{on $\Sigma$},\\
\partial_{\nu}u_2=0& \text{on $\Gamma$,}
\end{array}
\right. 
\end{equation}
with $\l_i\in\R$. Although mathematically it makes sense to consider $\l_i$ as a real parameter, its usual meaning is that $1/\l_i$ is the diffusion coefficient in $\O_i$, $\l_i$  being a positive parameter in such a case.

As a first step towards the study of (\ref{sisintro}), it is necessary  to analyze the eigenvalue problem 

\begin{equation}
\label{sisintroeigen}
 \left\{\begin{array}{ll}
-\D u_i=\l m_i(x)u_i & \text{in $\Omega_i$},\\
u_i>0 & \text{in $\Omega_i$},\\
\noalign{\smallskip}
\partial_{\nu}u_i=\gamma_i(u_2-u_1)& \text{on $\Sigma$},\\
\partial_{\nu}u_2=0& \text{on $\Gamma$,}
\end{array}
\right. 
\end{equation}
where $m_i\in L^\infty(\O_i)$, $m_i\not\equiv 0$ in $\O_i$. (\ref{sisintroeigen}) has been analyzed  in \cite{chinos} in the self-adjoint case $\g_1=\g_2$. For that, the authors used variational arguments to prove the existence of principal eigenvalue as well as its main properties.  The general case $\g_1\neq \g_2$ was studied in \cite{bcs} using a different argument. In \cite{bcs},  to study (\ref{sisintroeigen}), the authors first analyze the problem 
\begin{equation}
\label{sisintroeigenotro}
 \left\{\begin{array}{ll}
-\D u_i+c_i(x)u_i=\l u_i & \text{in $\Omega_i$},\\
\noalign{\smallskip}
\partial_{\nu}u_i=\gamma_i(u_2-u_1)& \text{on $\Sigma$},\\
\partial_{\nu}u_2=0& \text{on $\Gamma$,}
\end{array}
\right. 
\end{equation}
where $c_i\in L^\infty(\O_i)$. They prove the existence of a unique principal eigenvalue of  (\ref{sisintroeigenotro}), denoted by $\L_1(c_1,c_2)$. Hence, the study of (\ref{sisintroeigen}) is equivalent to find the zeros of the map
$$
\l\in\R\mapsto f(\l):=\L_1(-\l m_1,-\l m_2).
$$ 

The main goal of this paper is to study the following generalized eigenvalue problem:
\begin{equation}
\label{intronuestro}
 \left\{\begin{array}{ll}
-\D u_i=\l_i m_i(x)u_i & \text{in $\Omega_i$},\\
u_i>0 & \text{in $\Omega_i$},\\
\noalign{\smallskip}
\partial_{\nu}u_i=\gamma_i(u_2-u_1)& \text{on $\Sigma$},\\
\partial_{\nu}u_2=0& \text{on $\Gamma$.}
\end{array}
\right. 
\end{equation}
Motivated by \cite{bcs}, to study (\ref{intronuestro}) we analyze the zeros of the map
$$
(\l_1,\l_2) \in \R^2\mapsto F(\l_1,\l_2):= \L_1(-\l_1 m_1,-\l_2 m_2),
$$
that is, we analyze the set 
$$
{\cal C}:=\{(\l_1,\l_2)\in\R^2: F(\l_1,\l_2)=0\}.
$$
We show that $F$ is a regular function, concave and $F(0,0)=0$. Hence, for instance, fixed $\l_1$, there exist at most two values of $\l_2$ such that $F(\l_1,\l_2)=0$. Moreover, due to the concavity of $F$, it is well known that the set $\{(\l_1,\l_2)\in\R^2: F(\l_1,\l_2)\geq0\}$ is convex. In any case, the study of the set ${\cal C}$ depends strongly on the signs of $m_i$.   It is obvious that the case where the functions $m_i$ have a definite sign, for example they are positive, is a simpler case than the case where one or both of them change sign.

We summarize the main results.

\begin{figure}[h]
\centering
\includegraphics[width=.45\linewidth]{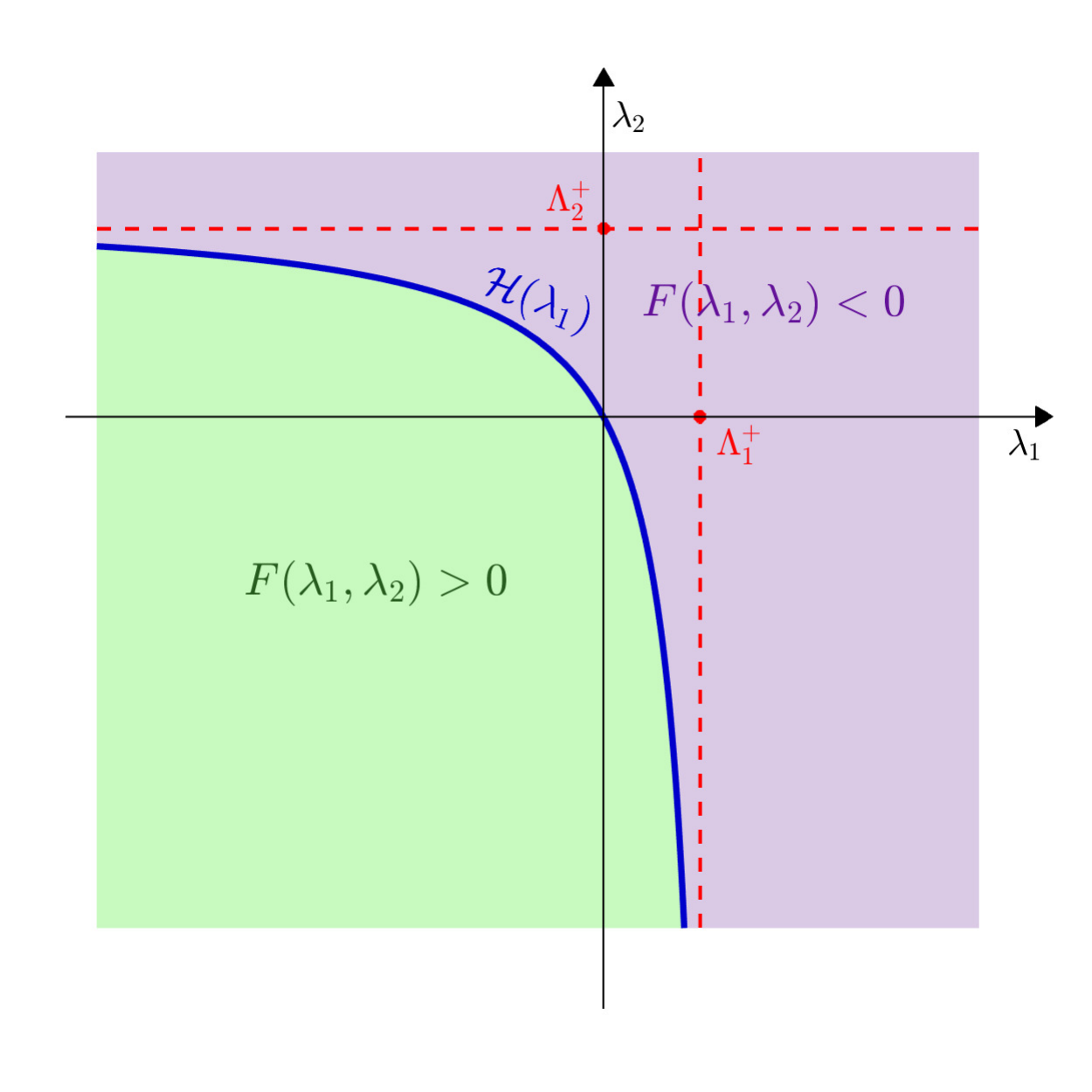}
\caption{Case $m_1$ and $m_2$ both non-negative and non-trivial verifying (\ref{00intro}): we have represented in the plane $\l_1-\l_2$ the curve $F(\l_1,\l_2)=0$, as well as, the regions where $F$ is negative and positive. In this case, ${\cal C}=\{(\l_1,\l_2): \l_2={\cal H}(\l_1)\}.$}
\label{figura:figure1}
\end{figure}

Our first result deals with the case both $m_i$  non-negative and non-trivial functions (see Figure~\ref{figura:figure1}). 
 \begin{theorem}
 \label{teoprinci1}
 Assume that $m_i\gneq 0$ in $\O_i$, $i=1,2$ and define
$$
M_i^0:=\Omega_i\setminus\overline{\{x\in \O_i:m_i(x)>0\}}
$$
and assume that  $\partial M_i^0$ is regular and 
\begin{equation}
\label{00intro}
\overline{M_i^0}\subseteq \O_i.
%dist(\partial M_i^0\cap \O_i,\partial \O_i)>0.
\end{equation}
Then, there exist positive values $\L_i^+$, $i=1,2$ such that:
 \begin{enumerate}
  \item Assume that $\l_1\geq \L_1^+$. Then, $F(\l_1,\l_2)<0$ for all $\l_2\in\R$.
 \item Assume that $\l_1<\L_1^+$. There exists a unique $\l_2:={\cal H}(\l_1)$ such that $F(\l_1,\l_2)=0$ and
 $$
F(\l_1,\l_2)<0\quad\mbox{for $\l_2>{\cal H}(\l_1)$},\quad  F(\l_1,\l_2)>0\quad\mbox{for $\l_2<{\cal H}(\l_1)$.}
 $$
 Moreover, the map $\l_1\mapsto {\cal H}(\l_1)$ is continuous, decreasing, ${\cal H}(0)=0$ and
 $$
 \lim_{\l_1\to -\infty}{\cal H}(\l_1)=\L_2^+,\qquad  \lim_{\l_1\to \L_1^+}{\cal H}(\l_1)=-\infty.
 $$
 \end{enumerate}
 \end{theorem}
 
The values $\L_1^+$ and $\L_2^+$ will be defined in Section 2.

\begin{figure}[h]
%\centering
%\begin{subfigure}{.5\textwidth}
\centering
\includegraphics[width=.45\linewidth]{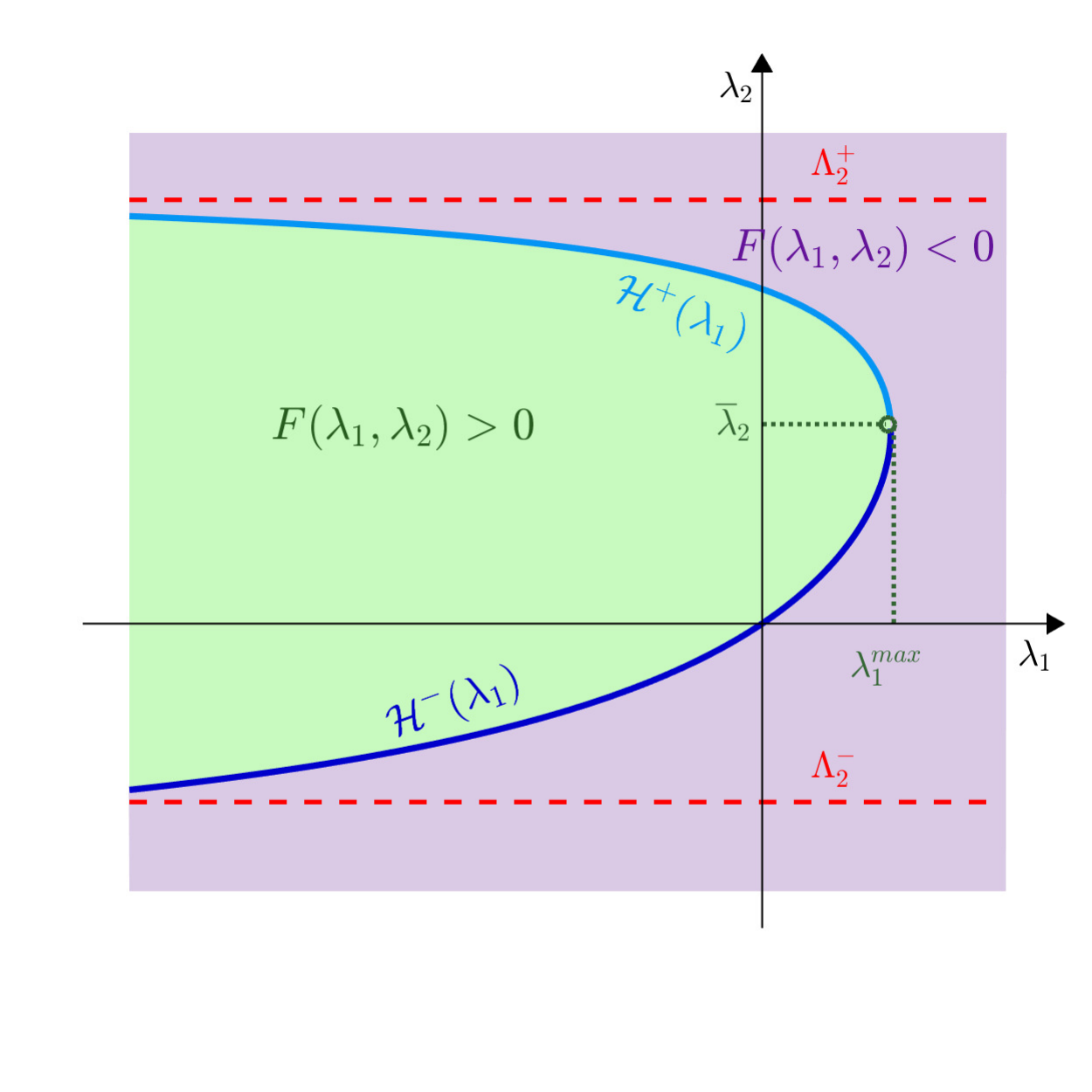}
\caption{Case $m_1\gneq  0$ and verifying (\ref{00intro}), $m_2$ changing sign and $\int_{\Omega_{2}}m_{2}<0$.}
\label{fig:sub1}
\end{figure}

In the next result we analyze the case $m_1$ non-negative and non-trivial and $m_2$ changing sign (see Figure \ref{fig:sub1}). 

\begin{theorem}
\label{teoprinci2} 
 Assume that $m_1\gneq 0$ in $\O_1$ and verifies (\ref{00intro}) and $m_2$ changes sign in $\O_2$. There exists $\l_1^{\max}\geq 0$ such that:
 \begin{enumerate}
 \item If $\l_1>\l_1^{\max}$, then $F(\l_1,\l_2)<0$ for all $\l_2\in\R$.
 \item If $\l_1=\l_1^{\max}$, then there exists a unique $\overline{\l}_2$ such that $F(\l_1^{\max},\overline{\l}_2)=0$  and  $F(\l_1^{\max},\l_2)<0$ for all $\l_2\in\R\setminus\{\overline{\l}_2\}$.
 \item  For all $\l_1<\l_1^{\max}$, there exist $\l_2^-={\cal H}^-(\l_1)<\l_2^+={\cal H}^+(\l_1)$ such that
 $$
 F(\l_1,\l_2^-)= F(\l_1,\l_2^+)=0,
 $$  
 and
 $$
 F(\l_1,\l_2)
 \left\{
 \begin{array}{ll}
 < 0 & \mbox{for $\l_2>{\cal H}^+(\l_1)$ or $\l_2<{\cal H}^-(\l_1)$},\\
 >0 & \mbox{for $\l_2\in ({\cal H}^-(\l_1),{\cal H}^+(\l_1))$.}
\end{array}
\right.
 $$
 Moreover, the map $\l_1\mapsto{\cal H}^+(\l_1)$ (resp. ${\cal H}^-(\l_1)$) is continuous, decreasing (resp. increasing) and
 $$
 \lim_{\l_1\to -\infty}{\cal H}^{\pm}(\l_1)=\L_2^\pm\quad\mbox{and}\quad  \lim_{\l_1\to \l_1^{\max}}{\cal H}^{\pm}(\l_1)=\overline{\l}_2.
 $$

\item Finally, 
\begin{enumerate}
\item If $\int_{\O_2} m_2<0$, then $ \l_1^{\max}>0$ and $\overline{\l}_2>0$.
\item If $\int_{\O_2} m_2>0$, then $ \l_1^{\max}>0$ and $\overline{\l}_2<0$.
\item If $\int_{\O_2} m_2=0$, then $\l_1^{\max}=\overline{\l}_2=0$.
\end{enumerate}

 \end{enumerate}

 \end{theorem} 

\begin{figure}
\centering
\includegraphics[width=.3\linewidth]{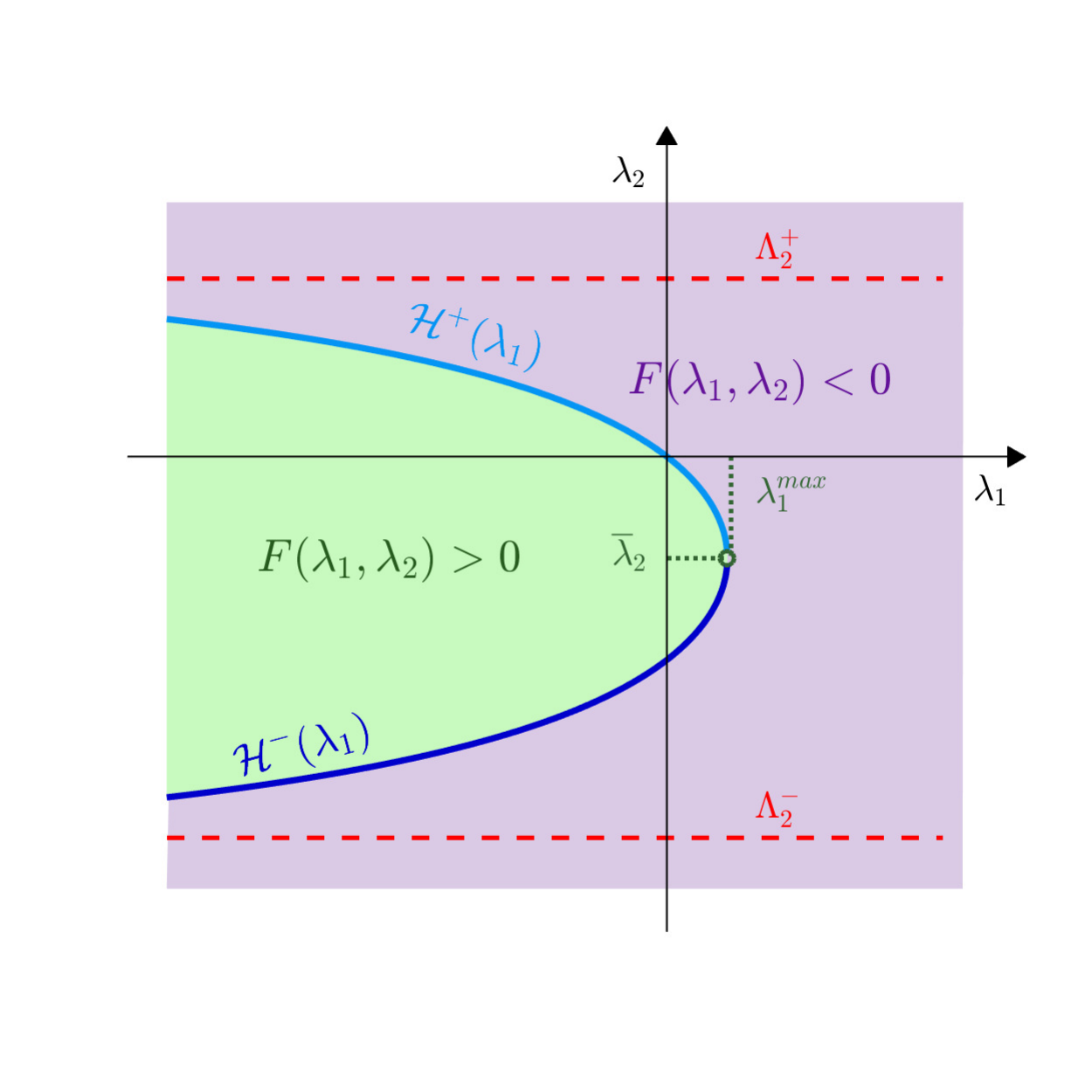}
\includegraphics[width=.3\linewidth]{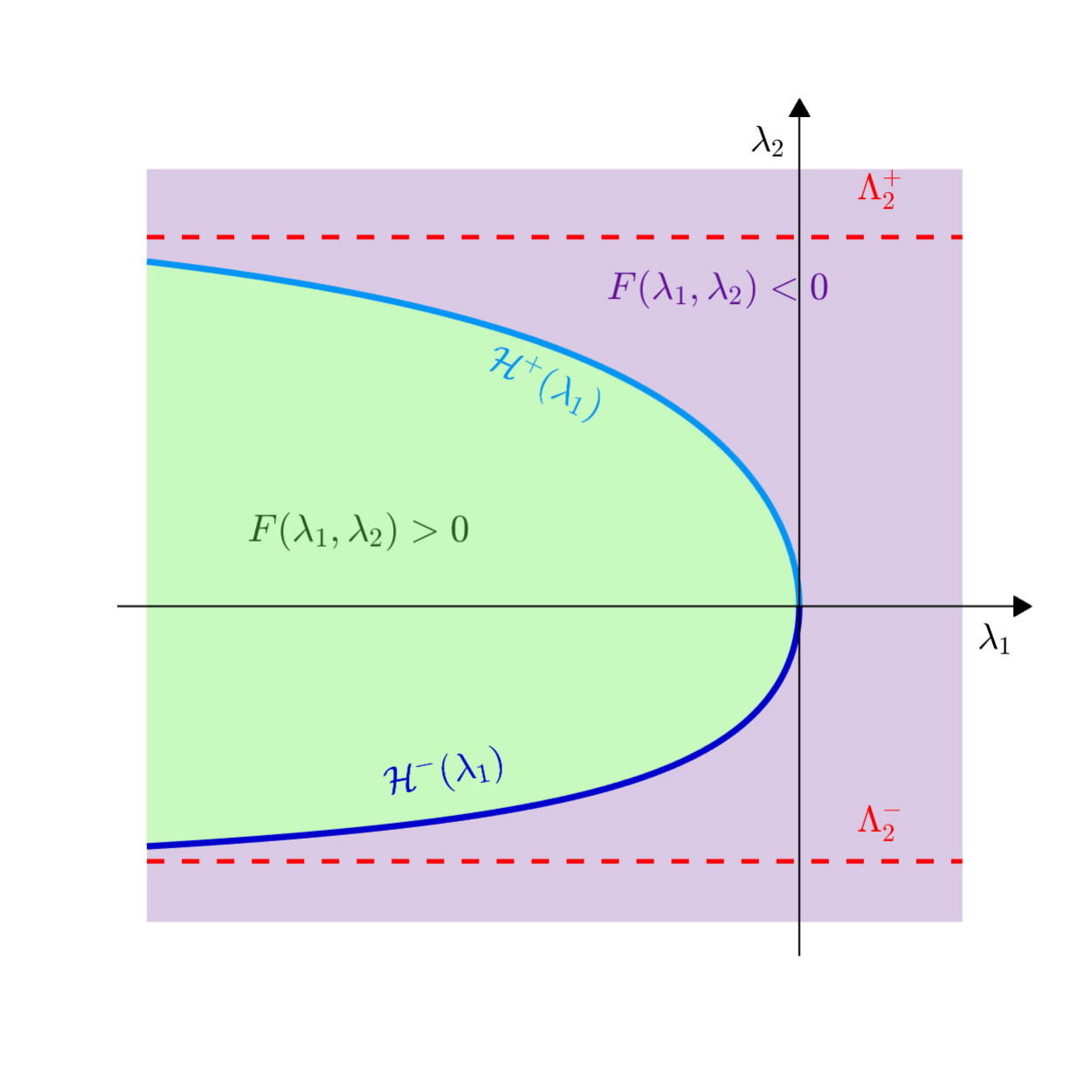}

\caption{Case $m_1\gneq  0$ and verifying (\ref{00intro}), $m_2$ changing sign and $\int_{\Omega_{2}}m_{2}>0$ (left) and $\int_{\Omega_{2}}m_{2}=0$ (right).}
\label{fig:sub2}
\end{figure}

\begin{remark}
\begin{enumerate}
\item In Figure \ref{fig:sub2} we have represented the cases $m_1\gneq 0$ in $\O_1$ verifying (\ref{00intro}), $m_2$ changes sign in $\O_2$,  $\int_{\O_2}m_2>0$ and  $\int_{\O_2}m_2=0$.
\item Of course, by symmetry, a similar result holds for $m_1$ changing sign in $\O_1$ and $m_2$ non-negative, non-trivial and verifying (\ref{00intro}).
\end{enumerate}
\end{remark}

Finally, we deal with the case of both $m_i$ changing sing.
 \begin{theorem}
 \label{teo3}
 Assume that $m_i$ changes sign in $\O_i$. Then, there exists a closed curve ${\cal C}\subset \R^2$, such that $F(\l_1,\l_2)=0$ if and only if $(\l_1,\l_2)\in {\cal C}$. Moreover,
 $$
\mbox{$F(\l_1,\l_2)>0$ if and only if $(\l_1,\l_2)\in int({\cal C})$,}
$$
and
$$
 \mbox{$F(\l_1,\l_2)<0$ if and only if $(\l_1,\l_2)\in Ext({\cal C})$.}
 $$
 \end{theorem} 
 
 \begin{remark}
The form and structure of ${\cal C}$ depends strongly on the sign of the integrals of $m_i$. In all the cases, $(0,0)\in {\cal C}$. 
 \end{remark}
 
\begin{figure}
\centering
\includegraphics[width=.25\linewidth]{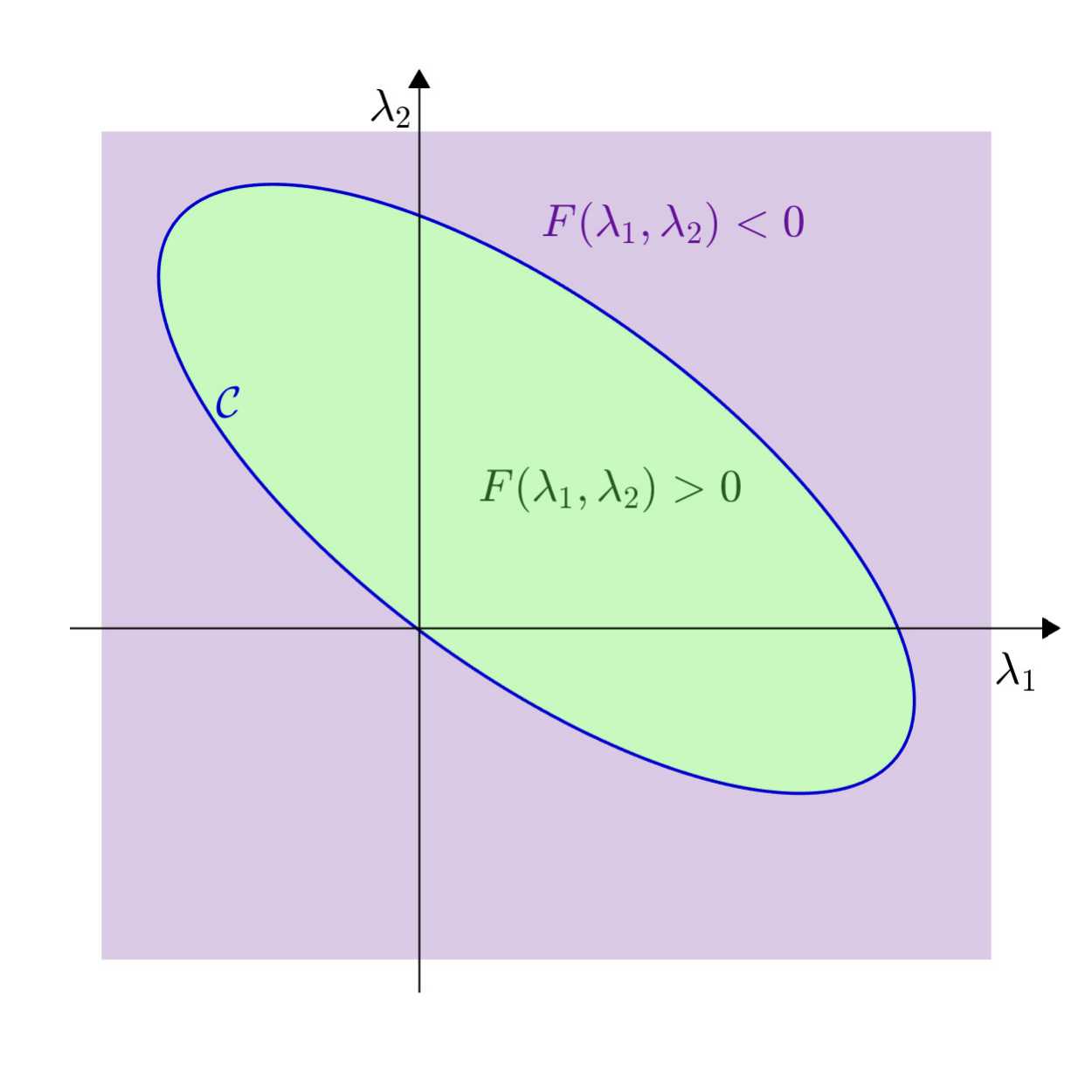}
\includegraphics[width=.25\linewidth]{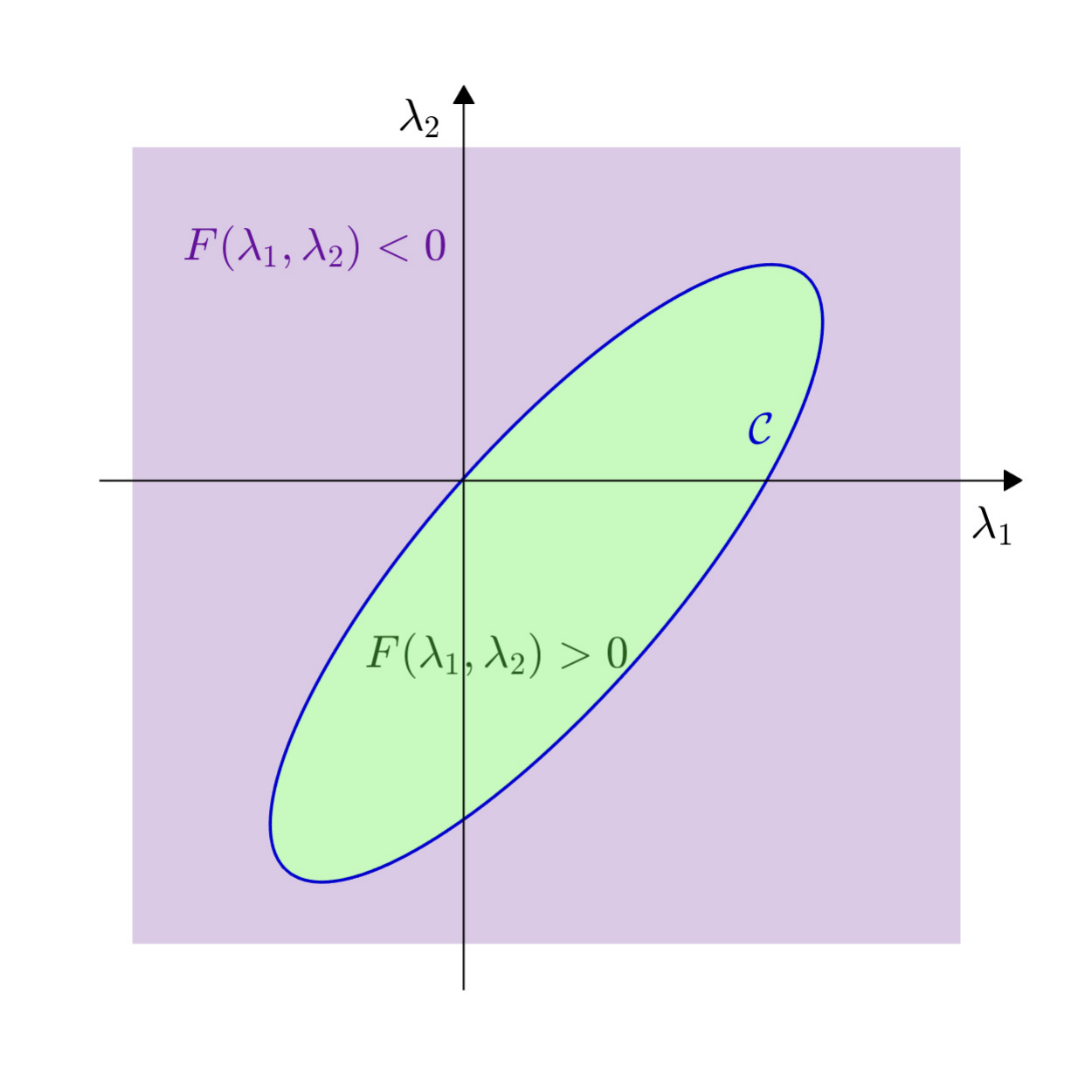}
\includegraphics[width=.25\linewidth]{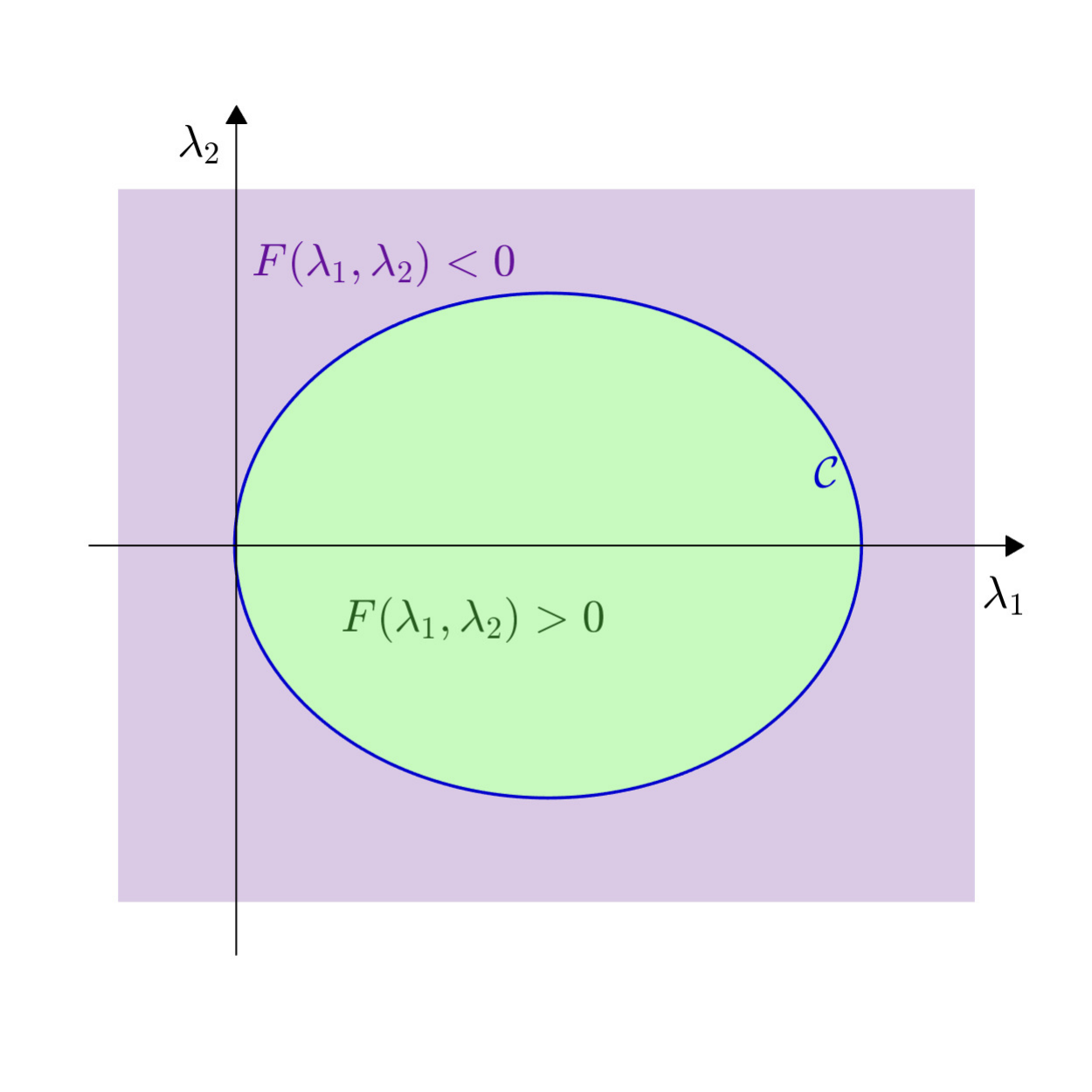}
\caption{Cases $m_i$ changing sign. $\int_{\Omega_1}m_1<0$ and: $\int_{\Omega_{2}}m_{2}<0$ (left), $\int_{\Omega_{2}}m_{2}>0$ (center) and $\int_{\Omega_{2}}m_{2}=0$ (right).}
\label{habi}
\end{figure}

 In the following result, we complete the above Theorem, see Figure \ref{habi}.

 \begin{theorem}
 \label{teoprinci3} 
 Assume that $m_i$ changes sign for $i=1,2$. There exist $\l_1^{\min}\leq 0\leq \l_1^{\max}$ such that
 \begin{enumerate}
  \item If $\l_1<\l_1^{\min}$ or $\l_1>\l_1^{\max}$, then $F(\l_1,\l_2)<0$ for all $\l_2\in\R$.
   \item If $\l_1=\l_1^{\max}$ (resp. $\l_1=\l_1^{\min}$) then there exists a unique $\overline{\l}_2$ (resp. $\underline{\l}_2$) such that $F(\l_1^{\max},\overline{\l}_2)=0$  (resp. $F(\l_1^{\min},\underline{\l}_2)=0$) and  $F(\l_1^{\max},\l_2)<0$ (resp. $F(\l_1^{\min},\l_2)<0$) for all $\lambda_2 \in\R\setminus\{\overline{\lambda}_2\}$(resp. $\lambda_2 \in\R\setminus\{\underline{\lambda}_2\}$).
 \item If $\l_1\in (\l_1^{\min}, \l_1^{\max})$ there exist unique $\l_2^-={\cal H}^-(\l_1)<\l_2^+={\cal H}^+(\l_1)$ such that $F(\l_1,\l_2^\pm)=0$. Moreover,  
  $$
 F(\l_1,\l_2)
 \left\{
 \begin{array}{ll}
 < 0 & \mbox{for $\l_2>{\cal H}^+(\l_1)$ or $\l_2<{\cal H}^-(\l_1)$},\\
 >0 & \mbox{for $\l_2\in ({\cal H}^-(\l_1),{\cal H}^+(\l_1))$.}
\end{array}
\right.
 $$
 Finally,
 $$
 \lim_{\l_1\to\l_1^{\min}}{\cal H}^\pm(\l_1)=\underline{\l}_2,\quad\mbox{and}\quad  \lim_{\l_1\to\l_1^{\max}}{\cal H}^\pm(\l_1)=\overline{\l}_2.
 $$

 \end{enumerate}
 \end{theorem} 

We apply these results to the nonlinear problem
\begin{equation}
\label{logintro}
 \left\{\begin{array}{ll}
-\D u_i=\l_i m_i(x)u_i-u_i^{p_i} & \text{in $\Omega_i$},\\
\noalign{\smallskip}
\partial_{\nu}u_i=\gamma_i(u_2-u_1)& \text{on $\Sigma$},\\
\partial_{\nu}u_2=0& \text{on $\Gamma$,}
\end{array}
\right. 
\end{equation}
with $\l_i\in\R$, $m_i\in L^\infty (\O_i)$ and $p_i>1$. Again, we consider the general case $\l_i\in\R$ and $m_i\in L^\infty (\O_i)$ but of course from a biological point of view only some cases are interesting. For example, when $m_1\equiv 1$ in $\O_1$, $\l_1$ represents the growth rate of the species $u_1$, which could be positive or negative. When $m_1$ changes sign, only the case $\l_1>0$  should be considered.

We prove (Theorem \ref{logis})  that (\ref{logintro})  possesses a positive solution if and only if 
$$
F(\l_1,\l_2)<0.
$$
Moreover, in such a case, the solution is the unique positive solution. Hence, we can give the following consequences:
\begin{enumerate}
\item Assume that $m_1$ and $m_2$ are non-negative and non-trivial functions. 
\begin{enumerate}
\item For $\l_1$ large ($\l_1>\L_1^+$), there exists a positive solution for all $\l_2\in\R$.
\item For $\l_1<\L_1^+$, there exists a value $\l_2={\cal H}(\l_1)$ such that  
(\ref{logintro}) possesses a positive solution for $\l_2>{\cal H}(\l_1)$. 
\end{enumerate}

In both cases, for $\l_1>0$ we have that there exists a positive solution for negative growth rate ($\l_2$) of $u_2$. In the case without interface, this is not possible, that is, even if the population has negative growth in one part of the domain, the interface effect makes it possible for the species to persist throughout the domain.

\item Assume that $m_1$ is non-negative and non-trivial and $m_2$ changes sign. Then, if $\l_1$ is large, then there exists positive solution for all $\l_2\in\R$. However, for $\l_1<\l_1^{\max}$, then there exists positive solution for $\l_2<{\cal H}^-(\l_1)$  or $\l_2>{\cal H}^+(\l_1)$.
\item Assume that $m_1$ and $m_2$ change sign. There exist $\l_1^{\min}<\l_1^{\max}$ such that for $\l_1>\l_1^{\max}$ or $\l_1<\l_1^{\min}$, (\ref{logintro}) possesses a positive solution for all $\l_2\in \R$. However, for $\l_1\in (\l_1^{\min},\l_1^{\max})$, there exist ${\cal H}^-(\l_1)<{\cal H}^+(\l_1)$ such that (\ref{logintro}) possesses a positive solution only for $\l_2<{\cal H}^-(\l_1)$  or $\l_2>{\cal H}^+(\l_1)$.
  
\end{enumerate}

An outline of the paper is: in Section 2 we include some preliminary results related to scalar eigenvalue problems. Section 3 is devoted to show some general properties of $F(\l_1,\l_2)$. The main results concerning to the eigenvalue problem (\ref{intronuestro}) are proved in Section 4. Finally, in Section 5, we analyze (\ref{logintro}).

\section{Preliminary results}
\subsection{Scalar eigenvalue problem}

In this section we recall some results concerning to scalar eigenvalue problems, see \cite{santijlg} for example.   

Here $G$ is a $C^{2,\alpha}$, $\alpha\in (0,1)$, domain of $\R^N$, $\partial G =
\Gamma_1 \cup \Gamma_2$, where $\Gamma_1$ and $\Gamma_2 $ are two disjoint open and closed subsets of 
$\partial G$ and $\nu$ is the outward unit normal vector field. For $c\in L^\infty(G)$, $h\in
C(\Gamma_1)$, $g\in C(\Gamma_2)$, we denote by $\sigma_1^G(-\D+c;N+h,N+g)$ the
principal eigenvalue of the problem
$$
\left\{
\begin{array}{lcc}
-\Delta \phi +c(x)\phi = \lambda \phi & \hbox{in} & G, \\[0.5pc]
\frac{\partial \phi }{\partial \nu}+h\phi = 0  & \hbox{on} & \Gamma_1, \\[0.5pc]
 \frac{\partial \phi }{\partial \nu}+g\phi = 0 & \hbox{on}
& \Gamma_2,
\end{array}
\right.
$$
and by $\sigma_1^G(-\D+c;N+h,D)$ that of the problem
$$
\left\{
\begin{array}{lcc}
-\Delta \phi  +c(x)\phi = \lambda \phi & \hbox{in} & G, \\
\frac{\partial \phi }{\partial \nu}+h\phi = 0  & \hbox{on} & \Gamma_1, \\[0.5pc]
\phi = 0 & \hbox{on} & \Gamma_2.
\end{array}
\right.
$$
We will quote some important properties of $\sigma_1^G(-\D+c;N+h,N+g)$ and
$\sigma_1^G(-\D+c;N+h,D)$. We denote the boundary operator
$$
{\cal B}(\phi)=\left\{
\begin{array}{lcc}
\partial_\nu  \phi +h\phi = 0  & \mbox{on  $\Gamma_1,$} \\
\partial_\nu \phi +g\phi = 0 & \mbox{on  $\Gamma_2$,}
\end{array}
\right.
\quad\mbox{or}\quad
{\cal B}(\phi)=\left\{
\begin{array}{lcc}
\partial_\nu  \phi +h\phi = 0  &  \mbox{on  $\Gamma_1,$} \\
\phi = 0 & \mbox{on  $\Gamma_2$.}
\end{array}
\right.
$$ 

\begin{proposition}
\label{propo1}
\begin{enumerate}
\item The map $c\in L^\infty(G)\mapsto \s_1^G(-\D+c;{\cal B}) $ is continuous and increasing.
\item It holds that
$$
\s_1^G(-\D+c;N+h, N+g)<\s_1^G(-\D+c; N+h, D).
$$
\item Assume that there exists a positive supersolution, that is, a positive function $\overline{u}\in W^{2,p}(G)$, $p>N$, such that
$$
-\D\overline{u}+c(x)\overline{u}\geq 0\quad\mbox{in $G$,}\quad {\cal B}(\overline{u})\geq 0\quad\mbox{on $\partial G$,}
$$
and some strict inequalities, then
$$
\s_1^G(-\D+c;{\cal B})>0.
$$
\end{enumerate}
\end{proposition}
Now, we define
$$
\mu(\l):=\sigma_1^G(-\D-\l c;{\cal B}),\quad\l\in\R.
$$
The main properties of $\m(\l)$ are stated in the next result.
\begin{proposition}
\label{propo2}
\begin{enumerate}
\item Assume that $c\not\equiv 0$ in $G$. Then, $\l\in\R\mapsto \m(\l)$ is regular and concave.
\item Assume that $c\gneq 0$ in $G$, define
\begin{equation}
\label{c0}
C_0:=G\setminus\overline{\{x\in G:c(x)>0\}},
\end{equation}
and assume that 
\begin{equation}
\label{c0bis}
\overline{C_0 }\subseteq G.
%dist(\partial C_0\cap G,\partial G)>0.
\end{equation}
The map $\l\mapsto \mu(\l)$ is decreasing and 
$$
\lim_{\l\to +\infty}\mu(\l)=-\infty,\qquad \lim_{\l\to -\infty}\mu(\l)=\sigma_1^{C_0}(-\D;D).
$$
\item Assume that $c$ changes sign, then 
$$
\lim_{\l\to \pm\infty}\mu(\l)=-\infty.
$$ 
Moreover, there exists $\l_0\in\R$ such that $\m'(\l_0)=0$, $\m'(\l)>0$ for $\l<\l_0$ and $\m'(\l)<0$ for $\l>\l_0$.
\end{enumerate}
\end{proposition}
We can describe  exactly the sign of the map $ \mu(\l)$.
\begin{corollary}
\label{importante}
\begin{enumerate}
\item Assume that $c\gneq 0$ and the set $C_0$ satisfies (\ref{c0bis}). Then, there exists a unique zero of the map $\m(\l)$, we denote it by $\l_1^+(G,c;{\cal B})$, and as consequence,
$$
\mu(\l)
\left\{
\begin{array}{ll}
>0 & \mbox{if $\l<\l_1^+(G,c;{\cal B})$},\\
=0 &  \mbox{if $\l=\l_1^+(G,c;{\cal B})$},\\
<0 & \mbox{if $\l>\l_1^+(G,c;{\cal B})$}.
\end{array}
\right.
$$
\item Asume that $c$ changes sign.
\begin{enumerate}
\item If $\m(\l_0)<0$, then $\m(\l)<0$ for all $\l\in\R$.
\item If $\m(\l_0)=0$, then $\l_0$ is the unique zero of the map $\m(\l)$.
\item If $\m(\l_0)>0$, then there exist two zeros of the map $\m(\l)$, we call them $\l_1^-(G,c;{\cal B})<\l_0<\l_1^+(G,c;{\cal B})$. As a consequence,
$$
\mu(\l)
\left\{
\begin{array}{ll}
>0 &  \mbox{if $\l\in (\l_1^-(G,c;{\cal B}),\l_1^+(G,c;{\cal B}))$}, \\
=0 &  \mbox{if $\l=\l_1^-(G,c;{\cal B})$ or $\l=\l_1^+(G,c;{\cal B})$},\\
<0 & \mbox{if $\l<\l_1^-(G,c;{\cal B})$ or $\l>\l_1^+(G,c;{\cal B})$}.
\end{array}
\right.
$$

 \end{enumerate}
 \end{enumerate}
\end{corollary}

\subsection{Interface eigenvalue problem}
First, we fix some notations that will be used throughout the paper. For convenience, we write ${\bf u}=(u_1,u_2)$ with $u_i$ defined in $\O_i$ and similarly ${\bf c}=(c_1,c_2)$. In order to simplify the notation we write the boundary conditions as
$$
\left\{\begin{array}{ll}
\partial_{\nu}u_i=\gamma_i(u_2-u_1)& \text{on $\Sigma$},\\
\partial_{\nu}u_2=0& \text{on $\Gamma$,}
\end{array}\right\}
\Longleftrightarrow
{\cal I}({\bf u})=0 \quad\mbox{on $\Sigma\cup \Gamma$.}
$$
We write
$$
{\cal I}(\bu)\succeq 0\quad\mbox{on $(\Sigma,\Gamma)$}\Longleftrightarrow
 \left\{\begin{array}{ll}
\partial_{\nu}u_1\geq\gamma_1(u_2-u_1)& \text{on $\Sigma$},\\
\partial_{\nu}u_2\leq\gamma_2(u_2-u_1)& \text{on $\Sigma$},\\
\partial_{\nu}u_2\geq 0& \text{on $\Gamma$.}
\end{array}
\right.
$$

We consider the Banach spaces
$$
\begin{array}{l}
L^p:=\left\{{\bf u}: u_i\in L^p(\O_i)\right\},\quad p\geq 1,\\
H^1:=\left\{{\bf u}: \mbox{$u_i\in H^1(\O_i)$}\right\},\\
W^{2,p}:=\left\{{\bf u}: u_i\in W^{2,p}(\O_i)\right\},\quad p\geq 1.\\
\end{array}
$$
The norm of a function $\bu$ is defined as the sum of the norms of $u_i$ in the respective spaces.

On the other hand, given $\bu=(u_1,u_2)$ we write $\bu\geq 0$ in $\O$ if $u_i\geq 0$ in $\O_i$ for $i=1,2$ and $\bu>0$ in $\O$ if both $u_i>0$ in $\O_i$ for $i=1,2$, and finally $\bu\neq 0$ in $\O$ if $u_i\neq 0$ in a subset of positive measure of $\O_i$ for some $i=1,2$.

Given $c_i\in L^\infty(\O_i)$, we denote by $\L_1({\bf c})=\L_1(c_1,c_2)$ the principal eigenvalue of (see \cite{chinos})  
\begin{equation}
\label{eigeninter}
 \left\{\begin{array}{ll}
-\D u_i+c_i(x)u_i=\l u_i& \text{in $\Omega_i$},\\
\noalign{\smallskip}
{\cal I}({\bf u})=0 & \mbox{on $\Sigma\cup \Gamma$.}
\end{array}
\right. 
\end{equation}

First, we recall some properties of $\L_1(c_1,c_2)$, see \cite{bcs}.  
\begin{definition}
\label{super}
Given $\overline\bu\geq 0$ in $\O$, $\overline\bu=(\overline u_1,\overline u_2)$, $\overline\bu\in W^{2,p}$, $p>N$, is a strict supersolution of $(-\D+{\bf c},{\cal I})$ if
$$
-\Delta\overline\bu +{\bf c}(x)\overline \bu\geq 0\quad \mbox{in $\O$,}\quad {\cal I}(\overline\bu)\succeq 0\quad \mbox{on $(\Sigma,\Gamma),$}
$$
and some of these inequalities are strict.
\end{definition}

\begin{proposition}
\label{monoconti}
\begin{enumerate}
\item Assume  that ${\bf c}\leq {\bf d}$ in $\O$. Then, $\L_1({\bf c})\leq \L_1({\bf d})$. Moreover, if ${\bf c}\neq {\bf d}$ in $\O$, $\L_1({\bf c})< \L_1({\bf d})$.
\item Assume that ${\bf c}_n\to {\bf c}$ in $L^{\infty}$, then $\L_1({\bf c}_n)\to \L_1({\bf c})$.
\item It holds that
$$
\L_1({\bf c})<\min\{\s_1^{\O_1}(-\D+c_1;N+\g_1),\s_1^{\O_2}(-\D+c_2;N+\g_2,N)\}.
$$
\item The map ${\bf c}\in L^\infty\mapsto \L_1({\bf c}) $ is concave.
\item $\L_1({\bf c})>0$ if and only if there exists a strict positive supersolution $\overline{\bf u}$ of $(-\D+{\bf c},{\cal I})$.
\end{enumerate}
\end{proposition}

%Moreover, given two changing sign functions $m_i\in L^\infty(\O_i)$, denote by $\Sigma(m_1,m_2)$ the principal eigenvalue of
%\begin{equation}
%\label{eigeninter3}
% \left\{\begin{array}{ll}
%-\D u_i=\l m_i(x)u_i& \text{in $\Omega_i$},\\
%\noalign{\smallskip}
%{\cal I}({\bf u})=0 & \mbox{on $\Sigma\cup \Gamma$.}
%\end{array}
%\right. 
%\end{equation}
%It is clear that if we define 
%$$
%F(\l):=\Lambda_1(-\l m_1,-\l m_2),
%$$
%then, 
%$$
%F(0)=F(\Sigma(m_1,m_2))=0.
%$$
%In fact, it is know, see \cite{bcs}, that $F(\l)$ is a concave function and $F(\l)>0$ for $\l\in (0,\Sigma(m_1,m_2))$ and $F(\l)<0$ for $\l\in (-\infty,0)\cap (\Sigma(m_1,m_2),\infty)$.

\section{Generalized interface principal eigenvalue: first properties}
The main goal in this paper is to analyze the eigenvalue problem
\begin{equation}
\label{eigeninterour}
 \left\{\begin{array}{ll}
-\D u_i=\l_i m_i(x)u_i& \text{in $\Omega_i$},\\
u_i>0 & \text{in $\Omega_i$},\\
\noalign{\smallskip}
{\cal I}({\bf u})=0 & \mbox{on $\Sigma\cup \Gamma$,}
\end{array}
\right. 
\end{equation}
where $m_i\in L^\infty(\O_i)$, $m_i\not\equiv 0$ in $\O_i$, $i=1,2$. It is obvious that
$$
\mbox{$(\l_1,\l_2)$ is an eigenvalue of (\ref{eigeninterour}) if and only if $\Lambda_1(-\l_1 m_1,-\l_2 m_2)=0$.}
$$
Hence, we define $F:\R^2\mapsto\R$ by
$$
F(\l_1,\l_2):=\Lambda_1(-\l_1 m_1,-\l_2 m_2).
$$
%The main goal of this paper is to study the zeros of $F$ as well as its sign.

The following result addresses the concavity of $\L_1(c_1,c_2)$ in each component.
\begin{proposition}
\label{concave1}
Fix $c_2\in L^\infty(\O_2)$. Then, the map $c_1\in L^\infty(\O_1)\mapsto \L_1(c_1,c_2)\in\R$ is concave.
\end{proposition}
\begin{proof}
Denote $G(c_1):=\L_1(c_1,c_2)$, take $c^i_1\in L^\infty(\O_1)$, $i=1,2$ and $t\in [0,1]$. Then, 
$$
G(tc_1^1+(1-t)c_1^2)=\L_1(tc_1^1+(1-t)c_1^2,c_2)=\L_1(tc_1^1+(1-t)c_1^2,tc_2+(1-t)c_2)=\L_1(t{\bf c}+(1-t){\bf d}),
$$
where ${\bf c}=(c_1^1,c_2)$ and ${\bf d}=(c_1^2,c_2)$.
Using now Proposition \ref{monoconti} 4., we get that
$$
\begin{array}{rl}
G(tc_1^1+(1-t)c_1^2)= & \L_1(t{\bf c}+(1-t){\bf d}) \\
\geq  &t\L_1({\bf c})+(1-t)\L_1({\bf d})\\
= & t\L_1(c_1^1,c_2)+(1-t)\L_1(c_1^2,c_2)\\
= & tG(c_1^1)+(1-t)G(c_1^2).
\end{array}
$$
This completes the proof.
\end{proof}
As a consequence, we deduce the concavity of the map $F(\l_1,\l_2)$.
\begin{corollary}
\label{concavo}
Fixed $\l_1\in\R$, $\l_2\mapsto F(\l_1,\l_2)$ is concave, and then, there exist at most two values of $\l_2$  such that $F(\l_1,\l_2)=0$. A similar result holds when we fix $\l_2$.
\end{corollary}
In order to simplify the notation, we denote (recall Corollary \ref{importante})
\begin{equation}
\label{masmenos}
\L_1^\pm:=\l_1^\pm(\O_1,m_1;N+\g_1),\qquad \L_2^\pm:=\l_1^\pm(\O_2,m_2;N+\g_2,N).
\end{equation}
Observe that if we denote by $\m(\l)=\s_1^{\O_1}(-\D-\l m_1;N+\g_1)$, then $\m(0)>0$. Hence, if $m_1$ changes sign the existence of $\L_1^-<0<\L_1^+$ is guaranteed by Corollary \ref{importante}. If $m_1\gneq 0$ in $\O_1$ then $\L_1^-=-\infty$.

The first result provides upper bounds of  $F(\l_1,\l_2)$.
\begin{lemma}
\label{cotas}
It holds:
\begin{equation}
\label{cota1}
F(\l_1,\l_2)< \min\{\s_1^{\O_1}(-\D-\l_1 m_1;N+\gamma_1),\s_1^{\O_2}(-\D-\l_2 m_2;N+\gamma_2,N)\},
\end{equation}
and
\begin{equation}
\label{cota2}
F(\l_1,\l_2)\leq \frac{\displaystyle -\l_1\int_{\O_1}m_1-\l_2\int_{\O_2}m_2+(\g_1+\g_2)|\Sigma |}{|\O_1|+|\O_2|}.
\end{equation}
 \end{lemma}
 \begin{proof}
 (\ref{cota1}) follows from Proposition \ref{monoconti} 3.
 
 Let $\varphi=(\v_1,\v_2)$ be a positive eigenfunction associated to $F(\l_1,\l_2)$. Observe that
 $$
 \begin{array}{l}
 -\D \v_i-\l_im_i(x)\v_i=F(\l_1,\l_2)\v_i\quad\mbox{in $\O_i$,}\\
{\cal I}(\varphi)=0 \quad \mbox{on $\Sigma\cup \Gamma$.}
 \end{array}
 $$
 Multiplying  by $1/\v_i$, integrating and adding the two resulting equations, we obtain
 $$
F(\l_1,\l_2)(|\O_1|+|\O_2|)=-\l_1\int_{\O_1}m_1-\l_2\int_{\O_2}m_2-\left(\int_{\O_1}\frac{|\n\v_1|^2}{\v_1^2}+ \int_{\O_2}\frac{|\n\v_2|^2}{\v_2^2}\right)
 $$
 $$
 +\int_\Sigma (\v_2-\v_1)\left(\frac{\g_2}{\v_2}-\frac{\g_1}{\v_1}\right).
 $$
 Observe that
 $$
 \int_\Sigma (\v_2-\v_1)\left(\frac{\g_2}{\v_2}-\frac{\g_1}{\v_1}\right)=(\g_1+\g_2)|\Sigma|-\int_\Sigma\frac{\g_1\v_2^2+\g_2\v_1^2}{\v_1\v_2},
 $$
 whence we conclude (\ref{cota2}).

 \end{proof}
 \begin{corollary}
\label{cotaraiz}
\begin{enumerate}
\item $F(\l_1,\l_2)<0$ for $\l_1\in (-\infty,\L_1^-]\cup [\L_1^+,+\infty)$ or $\l_2\in (-\infty,\L_2^-]\cup [\L_2^+,+\infty)$.
\item Assume  that $F(\l_1,\l_2)=0$. Then,
$$
\L_1^-<\l_1<\L_1^+\quad\mbox{and}\quad \L_2^-<\l_2<\L_2^+.
$$
\end{enumerate}
\end{corollary}
\begin{proof}
\begin{enumerate}
\item Observe that if $\l_1\in (-\infty,\L_1^-]\cup [\L_1^+,+\infty)$, then, by Corollary \ref{importante}, we get  $\s_1^{\O_1}(-\D-\l_1 m_1;N+\gamma_1)\leq 0$. Hence, by (\ref{cota1}) we obtain that
$$
F(\l_1,\l_2)<0.
$$
\item Since $F(\l_1,\l_2)=0$, by (\ref{cota1}), we have that 
$$
0<\min\{\s_1^{\O_1}(-\D-\l_1 m_1;N+\gamma_1),\s_1^{\O_2}(-\D-\l_2 m_2;N+\gamma_2,N)\},
$$
and then $\s_1^{\O_1}(-\D-\l_1 m_1;N+\gamma_1) >0$ and $\s_1^{\O_2}(-\D-\l_2 m_2;N+\gamma_2,N)>0$, and hence the result concludes by Corollary \ref{importante}.
\end{enumerate}
\end{proof}

The following result will be very useful. 
%First, we introduce some notation. Assume  that ${\bf m}\geq 0$ in $\O$, define
%$$
%M_i^0:=\Omega_i\setminus\overline{\{x\in \O_i:m_i(x)>0\}}
%$$
%and assume that  $\partial M_i^0$ is regular and 
%\begin{equation}
%\label{00}
%\overline{M_i^0}\subseteq \O_i.
%dist(\partial M_i^0\cap \O_i,\partial \O_i)>0.
%\end{equation}
 \begin{proposition}
 \label{convergence} 
 Assume that $m_i\gneq 0$ in $\O_i$ and that the set $M_2^0$ verifies (\ref{00intro}). Take two sequences $\{a_n\}$ and $\{b_n\}$ such that
 $$
 a_n\to a_*\in (-\infty,\infty),\quad b_n\to -\infty  \quad\mbox{as $n\to +\infty$.}
 $$
 Then, at least for a subsequence, 
 $$
\lim_{n\to\infty}F(a_n,b_n) = \min\{\s_1^{\O_1}(-\D-a_*m_1;N+\g_1),\s_1^{M_2^0}(-\D;D)\}.
 $$
 \end{proposition}
 \begin{proof}
 Observe that
 $$
 F(a_n,b_n)<\min\{\s_1^{\O_1}(-\D-a_n m_1;N+\gamma_1),\s_1^{\O_2}(-\D-b_n m_2;N+\gamma_2,N)\} .
 $$
 By continuity,
 $$
  \s_1^{\O_1}(-\D-a_n m_1;N+\gamma_1)\to  \s_1^{\O_1}(-\D-a_* m_1;N+\gamma_1),
$$
and using Proposition \ref{propo2} 2., we get
$$
 \s_1^{\O_2}(-\D-b_n m_2;N+\gamma_2,N)\to \s_1^{M_2^0}(-\D;D).
 $$
 Hence, $F(a_n,b_n)$ is bounded.

 Assume that 
 \begin{equation}
\label{supo}
\s_0 :=\s_1^{M_2^0}(-\D;D)<\s_1^{\O_1}(-\D-a_* m_1;N+\gamma_1).
\end{equation}
Consequently, we conclude that, for a subsequence,  
$$
F(a_n,b_n)\to F_0\leq \s_0=\s_1^{M_2^0}(-\D;D)<\s_1^{\O_1}(-\D-a_* m_1;N+\gamma_1)<\infty \quad\mbox{as $n\to \infty$.}
$$ 
Without loss of generality, we consider ${\bf \v}_n=(\v_{1n},\v_{2n})$ a positive eigenfunction associated to $F(a_n,b_n)$ such that $\|{\bf \v}_n\|_2=1$.  Then,
$$
\int_\O |\nabla \v_n|^2-a_n\int_{\O_1} m_1\v_{1n}^2-b_n\int_{\O_2} m_2\v_{2n}^2+\int_\Sigma (\g_1\v_{1n}^2+\g_2\v_{2n}^2)-(\g_1+\g_2)\int_{\Sigma}\v_{1n}\v_{2n}=F(a_n,b_n)\leq C,
$$
where we have denoted
$$
\int_\O |\nabla \v_n|^2=\sum_{i=1}^2\int_{\O_i}|\nabla \v_{in}|^2.
$$
Since $b_n<0$, ${\bf m}\geq 0$ in $\O$ and $a_n\to a_*\in (-\infty,\infty)$, we get that
\begin{equation}
\label{aire}
\int_\O |\nabla \v_n|^2-(\g_1+\g_2)\int_{\Sigma}\v_{1n}\v_{2n}\leq C.
\end{equation}
Using now the inequalities
\begin{equation}
\label{desiuno}
\int_\Sigma u_1u_2\leq \frac{1}{2}\left(\int_\Sigma u_1^2+\int_\Sigma u_2^2\right),
\end{equation}
and that for any $\varepsilon>0$ there exists $C(\varepsilon)>0$ such that
\begin{equation}
\label{desidos}
\int_\Sigma v^2\leq \varepsilon\int_{\O_i}|\nabla v|^2  +C(\varepsilon)\int_{\O_i} v^2\quad\forall v\in H^1(\O_i),
\end{equation}
(see for instance Lemma 1 in \cite{ab99}) and $\|{\bf \v}_n\|_2=1$, wet get that
$$
\int_{\Sigma}\v_{1n}\v_{2n}\leq \frac{1}{2}\left(\int_\Sigma\v_{1n}^2+\int_\Sigma\v_{2n}^2\right)\leq \frac{1}{2}\left(\varepsilon\left(\int_{\O_1}|\nabla\v_{1n}^2|+\int_{\O_2}|\nabla\v_{2n}^2|\right)+C(\varepsilon)\right),
$$
and then from (\ref{aire}) we get
$$
\int_\O |\nabla \v_n|^2\leq C.
$$
Hence,
$$
\|\v_n\|_{H^1}\leq C_0.
$$
Thus, 
$$
\v_n\rightharpoonup\v_\infty=(\v_{1\infty},\v_{2\infty})\geq 0 \quad\mbox{in $H^1$},\quad \v_n\to\v_\infty \quad\mbox{in $L^2$ and $L^2(\Sigma)$ with $\|\v_\infty\|_2=1$.}
$$
By definition of $F(a_n,b_n)$ we have that
\begin{equation}
\label{debil}
\begin{array}{c}
\displaystyle\sum_{i=1}^2\left(\int_{\O_i}\nabla \v_{in}\cdot\nabla v_i -a_n\int_{\O_1} m_1\v_{1n}v_1-b_n\int_{\O_2} m_2\v_{2n}v_2\right)+
\\
\displaystyle +\int_{\Sigma}(\v_{2n}-\v_{1n})(\g_2v_2-\g_1 v_1)=F(a_n,b_n)\left(\int_{\O_1}\v_{1n}v_1+\int_{\O_2}\v_{2n}v_2 \right),\quad \forall v_i\in H^1(\O_i).
\end{array}
\end{equation}
First, we prove that
\begin{equation}
\label{ho1}
\v_{2\infty}\in H_0^1(M_2^0).
\end{equation}
Since 
$$
H_0^1(M_2^0)=\{u\in H^1(\O_2):\mbox{$u=0$ in $\O_2\setminus M_2^0$} \},
$$
we claim that $\v_{2\infty}=0$ in $\O_2\setminus M_2^0$.

By contradiction, assume that $\v_{2\infty}>0$ in $D$, for some $D\subset \O_2\setminus M_2^0$ and take $v_1=0$ in $\O_1$ and $v_2\in C_0^\infty(D)$, $v_2>0$ in $D$. Then, by (\ref{debil})
\begin{equation}
\label{sol}
-\int_D\D v_2\v_{2n}-b_n\int_D m_2(x)\v_{2n}v_2=F(a_n,b_n)\int_{D}\v_{2n}v_2.
\end{equation}
If $\v_{2\infty}>0$ in $D$, then $-b_n\int_D m_2(x)\v_{2n}v_2\to \infty$ as $b_n\to-\infty$, a contradiction with (\ref{sol}).
Hence, we conclude that $\v_{2\infty}=0$ in $D$. This implies (\ref{ho1}).
%$$
%\v_{2\infty}\in H_0^1(M_2^0).
%$$

Taking $v_1\in H^1(\O_1)$ and $v_2=0$ in (\ref{debil}), taking limit, we get
$$
\int_{\O_1}\nabla \v_{1n}\cdot\nabla v_1-a_n\int_{\O_1}m_1 \v_{1n}v_1+\displaystyle \int_{\Sigma}(\v_{2n}-\v_{1n})(-\g_1 v_1)
=F(a_n,b_n)\int_{\O_1}\v_{1n}v_1,
$$
then passing to the limit, taking into account (\ref{ho1}) in the boundary integral, 
$$
\int_{\O_1}\nabla \v_{1\infty}\cdot\nabla v_1-a_*\int_{\O_1}m_1 \v_{1\infty}v_1+\g_1\int_{\Sigma}\v_{1\infty} v_1=F_0\int_{\O_1} \v_{1\infty}v_1.
$$
Hence, if $\| \v_{1\infty}\|_2\neq 0$, then $F_0=\s_1^{\O_1}(-\D-a_*m_1;N+\g_1)$, an absurdum due to $F_0\leq\s_0$ and (\ref{supo}). Then,  $\| \v_{1\infty}\|_2=0$. Hence, 
$$
\| \v_{2\infty}\|_2= 1.
$$ 
Then, take $v_1=0$ and $v_2\in H_0^1(M_2^0)$  in (\ref{debil}), we obtain 
$$
\int_{M_0^2}\nabla \v_{2\infty}\cdot\nabla v_2=F_0\int_{M_0^2}\v_{2\infty}v_2,
$$
which yields that $F_0=\s_1^{M_2^0}(-\D;D)=\s_0$.

A similar reasoning can be carried out when $\s_1^{\O_1}(-\D-a_* m_1;N+\gamma_1)<\s_1^{M_2^0}(-\D;D)$.

This finishes the proof.

 \end{proof}

% Fix $\l_2\in\R$ and define
% $$
% f_1(\l_1):=F(\l_1,\l_2).
% $$
% \begin{proposition}
% \label{propiedades}
% \begin{enumerate}
% \item Assume that $m_1\geq 0$, $m_1\neq 0$ and $m_1=0$ in $\O_{01}\subset\O_1$ such that
% $$
% dist(\partial \O_{01}\cap\O_1,\Sigma)>0.
% $$
% Then, $\l_1\mapsto f_1(\l_1)$ is decreasing and
% $$
% \lim_{\l_1\to -\infty}f_1(\l_1)=\s_1^0,\quad\mbox{and}\quad  \lim_{\l_1\to +\infty}f_1(\l_1)=-\infty.
% $$
% \item Assume that $m_1$ changes sign. Then, 
%\begin{enumerate}
%\item $\l_1\mapsto f_1(\l_1)$ is concave.
%\item
%$$
%\lim_{\l_1\to \pm\infty}f_1(\l_1)=-\infty.
%$$
%\item Assume that $\l_2=0$, then $f_1(0)=0$. Moreover,
%\begin{enumerate}
%\item If $\int_{\O_1}m_1<0$, then there exists $\l_1^+>0$ such that
%$$
%f_1(\l_1)
%\left\{
%\begin{array}{ll}
%0 & \mbox{if $\l_1=\l_1^+$,}\\
%>0 & \mbox{if $\l_1\in(0,\l_1^+)$,}\\
%<0 & \mbox{if $\l_1\in(-\infty,0)\cup(\l_1^+,+\infty)$.}
%\end{array}
%\right.
%$$
%\item If $\int_{\O_1}m_1>0$, then there exists $\l_1^-<0$ such that
%$$
%f_1(\l_1)
%\left\{
%\begin{array}{ll}
%0 & \mbox{if $\l_1=\l_1^-$,}\\
%>0 & \mbox{if $\l_1\in(\l_1^-,0)$,}\\
%<0 & \mbox{if $\l_1\in(-\infty,\l_1^-)\cup(0,+\infty)$.}
%\end{array}
%\right.
%$$
%\item If $\int_{\O_1}m_1=0$, then $f_1(\l_1)>0$ for all $\l_1\neq 0$.
%\end{enumerate}
%\end{enumerate} 
% \end{enumerate}
% \end{proposition}

\section{Proofs of the main results}
The main idea of the proof can be summarized as follows. Instead of looking for solutions of $F(\l_1,\l_2)=0$ in the general form $(\l_1,\l_2)$, we look for solutions in the particular form $\l_2=\mu \l_1$, for all $\m\in\R$.

Hence, the following map plays an essential role in our study.  Given $\m\in \R$, we define
$$
f_\m(\l_1):=\L_1(-\l_1 m_1,-\l_1 \m m_2)=F(\l_1,\l_1\m).
$$
In the following result we state that, for $\l_1\neq 0$, it is equivalent to solve $F(\l_1,\l_2)=0$ to $f_\m(\l_1)=0$. Specifically, we have:
 
\begin{proposition} 
\label{equi}
Assume that $F(\l_1^0,\l_2^0)=0$ and $\l_1^0\neq 0$, then $f_{\m_0}(\l_1^0)=0$ for $\m_0=\l_2^0/\l_1^0$. 

Conversely, if $f_{\m_0}(\l_1^0)=0$ then 
$F(\l_1^0,\l_2^0)=0$ for $\l_2^0=\m_0\l_1^0$.
 \end{proposition}

In what follows, we explore the particular case $\l_1=0$.
\begin{proposition}
\label{caseparticular}
Assume that $\l_1=0$ and denote
$$
g(\l_2):=F(0,\l_2)=\L_1(0,-\l_2m_2).
$$
The map $\l_2\mapsto g(\l_2)$ is regular, concave, $g(0)=0$ and 
\begin{equation}
\label{derivadag}
sign(g'(0))=sign\left(-\int_{\O_2}m_2\right).
\end{equation}
Moreover,
\begin{enumerate}
\item If $m_2\gneq 0$ in $\O_2$ and $M_2^0$ verifies (\ref{00intro}), then $\l_2\mapsto g(\l_2)$ is decreasing and
$$
\lim_{\l_2\to +\infty}g(\l_2)=-\infty\quad\mbox{and}\quad \lim_{\l_2\to -\infty}g(\l_2)=\min\{\s_1^{\O_1}(-\D; N+\g_1),\s_1^{M_2^0}(-\D; D)\}.
$$
In this case, $g(\l_2)>0$ for $\l_2<0$ and $g(\l_2)<0$ if $\l_2>0$.
\item If $m_2$ changes sign in $\O_2$, then
$$
\lim_{\l_2\to \pm\infty}g(\l_2)=-\infty.
$$
Moreover, 
\begin{enumerate}
\item If $\int_{\O_2}m_2=0$, then  $g'(0)=0$ and $\l_2=0$ is the unique root of $g(\l_2)=0$. As a consequence,
$g(\l_2)<0$ for $\l_2\neq 0$.
\item If $\int_{\O_2}m_2<0$, then  $g'(0)>0$ and there exists $\l_2^+>0$ such that $g(\l_2^+)=0$. In this case,
$$
g(\l_2)
\left\{
\begin{array}{ll}
>0 &\mbox{if $\l_2\in (0,\l_2^+)$,}\\
<0  &\mbox{if $\l_2<0$ or  $\l>\l_2^+$.}
\end{array}
\right.
$$
\item If $\int_{\O_2}m_2>0$, then  $g'(0)<0$ and there exists $\l_2^-<0$ such that $g(\l_2^-)=0$. Hence,
$$
g(\l_2)
\left\{
\begin{array}{ll}
>0 &\mbox{if $\l_2\in (\l_2^-,0)$,}\\
<0  &\mbox{if $\l_2<\l_2^-$ or  $\l_2>0$.}
\end{array}
\right.
$$
\end{enumerate}
\end{enumerate} 
\end{proposition}
\begin{proof}
To begin with, the regularity of $g$ follows by the regularity of the function $F$. On the other hand,  by Proposition \ref{concave1} follows that $g(\l)$ is concave. It is obvious that $g(0)=F(0,0)=\L_1(0,0)=0$.  On the other hand, taking  ${\bf m}=(0,m_2)$ in Proposition 3.17 in \cite{bcs}, we conclude (\ref{derivadag}).

Finally, observe that by  (\ref{cota1}) we have
\begin{equation}
\label{moco}
g(\l_2)<\s_1^{\O_2}(-\D-\l_2 m_2;N+\g_2,N),
\end{equation}
whence we deduce that $\displaystyle\lim_{\l_2\to +\infty}g(\l_2)=-\infty$ from Proposition \ref{propo2} 2.  and 3.
\begin{enumerate}
\item Assume that $m_2\gneq 0$ in $\O_2$. In this case, $g$ is decreasing. Moreover, by Proposition \ref{convergence}, taking $a_n=0$, we conclude that
$$
\lim_{\l_2\to -\infty}g(\l_2)=\min\{\s_1^{\O_1}(-\D; N+\g_1),\s_1^{M_2^0}(-\D; D)\}.
$$
\item Assume that $m_2$ changes sign. Then, using (\ref{moco}) and  Proposition \ref{propo2} 3., we deduce that
$$
\lim_{\l_2\to -\infty}g(\l_2)=-\infty.
$$
%Denote by $\v_{\l_2}$ the positive eigenfunction associated to $g(\l_2)$. Then, 
%\begin{equation}
%\label{estopa}
%-\D\v_{\l_2}-\l{\bf m}\v_{\l_2}=g(\l_2)\v_{\l_2}\quad\mbox{in $\O$,}\quad {\cal I}(\v_{\l_2})=0\quad\mbox{on $(\Sigma,\Gamma)$,}
%\end{equation}
%where 
%$$
%{\bf m}=(0,m_2).
%$$
%
%After, differentiating with respect to $\l_2$, we get
%\begin{equation}
%\label{derivada}
%(-\D-\l_2{\bf m}-g(\l_2))\v'_{\l_2}=(g'(\l_2)+{\bf m})\v_{\l_2}\quad\mbox{in $\O$,}\quad {\cal I}(\v'_{\l_2})=0\quad\mbox{on $(\Sigma,\Gamma)$.}
%\end{equation}
%We need to calculate the adjoint of the eigenvalue problem (\ref{estopa}). It is not difficult to prove that 
%$$
%-\D\v^*_{\l_2}-\l{\bf m}\v^*_{\l_2}=g(\l_2)\v^*_{\l_2}\quad\mbox{in $\O$,}\quad {\cal I^*}(\v^*_{\l_2})=0\quad\mbox{on $(\Sigma,\Gamma)$,}
%$$
%where ${\cal I^*}(\bu)=0$ denotes
%$$
%\left\{
%\begin{array}{ll}
%\partial_\nu u_i=\g_2 u_2-\gamma_1u_1  & \mbox{on $\Sigma$,}\\
%\partial_\nu u_2=0 & \mbox{on $\Gamma$.}
%\end{array}
%\right.
%$$
%Then, multiplying (\ref{derivada}) by $\v_{\l_2}^*$ and integrating yields to
%$$
%g'(\l_2)=-\displaystyle\frac{\displaystyle\int_\O {\bf m}\v_\l\v_{\l_2^*}}{\displaystyle\int_\O \v_{\l_2}\v^*_{\l_2}}.
%$$
%Taking into account that $g(0)=0$ and that $\v_0=(1,1)$,  ${\bf m}=(0,m_2)$, and $\v_0^*=(1,\g_1/\g_2)$, it follows that
%$$
%g'(0)=-\frac{\g_1}{\g_2|\O_1|+\g_1|\O_2|}\int_{\O_2}m_2.
%$$
Now,  from the sign of $g'(0)$ in (\ref{derivadag}), we conclude the result.
\end{enumerate}
\end{proof}

In the next result, we study in detail the map $\l_1\mapsto f_\mu(\l_1)$.

\begin{proposition}
\label{casegeneral}
Fix $\m\in\R$. Then, $\l_1\mapsto f_\m(\l_1)$ is regular, concave, $f_\m(0)=0$ and 
\begin{equation}
\label{derivada}
sign(f_\m'(0))=-sign\left(\g_2\int_{\O_1}m_1+\m\g_1\int_{\O_2}m_2\right).
\end{equation}
\begin{enumerate}
\item If $m_1\gneq 0$ in $\O_1$, then
$$
\lim_{\l_1\to +\infty}f_\m(\l_1)=-\infty
$$
\item If $m_1$ or $m_2$ changes sign, then 
$$
\lim_{\l_1\to\pm \infty}f_\m(\l_1)=-\infty.
$$
\end{enumerate}
\end{proposition}
 \begin{proof}
 It is clear that $f_\mu(0)=0$. The regularity follows by the regularity of $F$, the concavity of $f_\mu(\l_1)$ follows by Proposition \ref{monoconti} 4., and  (\ref{derivada}) follows taking  ${\bf m}=(m_1,\m m_2)$ in Proposition 3.17 in \cite{bcs}.
 
 On the other hand, by (\ref{cota1}) we get 
 $$
 f_\m(\l_1)<\min\{\s_1^{\O_1}(-\D-\l_1 m_1;N+\g_1),\s_1^{\O_2}(-\D-\l_1\mu m_2;N+\g_2,N)\},
 $$ 
 and then $\displaystyle\lim_{\l_1\to +\infty}f_\m(\l_1)=-\infty$, and if $m_1$ or $m_2$ changes sign, $\displaystyle\lim_{\l_1\to -\infty}f_\m(\l_1)=-\infty.$
\end{proof}

For $\int_{\O_2}m_2\neq 0$, we define
\begin{equation}
\label{m1}
 \m^*:=-\frac{\displaystyle\g_2\int_{\O_1}m_1}{\displaystyle\g_1\int_{\O_2}m_2},
\end{equation}
 in such a way that $f'_{\m^*}(0)=0$.

\subsection{Case $m_i\gneq 0$ in $\O_i$, $i=1,2$.}
Observe that in this case
$$
\m^*<0.
$$ 
\begin{proposition}
\label{casegeneraligualsigno}
Assume that $m_i\gneq 0$ in $\O_i$ and $M_i^0$ verify (\ref{00intro}) for $i=1,2$. 
 \begin{enumerate}
 \item If $\m\geq 0 $, the unique zero of $f_\m(\l_1)$ is $\l_1=0$.
 \item If $\m<0$, $\m\neq \m^*$, there exists an unique $\l_1=h_1(\m)\neq 0$ such that $f_\m(\l_1)=0$. Moreover, 
 $$
 h_1(\m)
 \left\{
 \begin{array}{ll}
 <0 & \mbox{if $\m>\m^*$ },\\
 =0 & \mbox{if $\m=\m^*$,}\\
 >0 & \mbox{if $\m<\m^*$.}
 \end{array}  
 \right.
 $$
 \item The map $\m\in (-\infty, 0)\mapsto h_1(\m)$ is continuous and decreasing. Moreover,
 $$
 \lim_{\m\uparrow 0}h_1(\m)=-\infty,\quad  \lim_{\m\to -\infty}h_1(\m)=\L_1^+.
 $$ 
% Finally, given $c\in \R$ the equation
% $$
% \l_1(m)=c
% $$ 
% has at most two solutions.
 \end{enumerate}
 \end{proposition}
 \begin{proof}
  \begin{enumerate}
 \item Assume that $\m\geq 0$. Then, since $\l_1\mapsto f_\m(\l_1)$ is decreasing and $f_\m(0)=0$, the result follows.
 \item Assume that $\m<0$.  Recall that $\l_1\mapsto f_\m(\l_1)$ is concave and $f_\m(0)=0$. If $\m>\m^*$ then $f'_\m(0)<0$, and hence there exists a unique $h_1(\m)<0$ such that $f_\m(h_1(\m))=0$. Similarly, when $\m<\m^*$ there exists a unique $h_1(\m)>0$ such that $f_\m(h_1(\m))=0$. 
 
\item We will show that  $\m\mapsto h_1(\m) $ is decreasing. Take now $\m_1<\m_2<0$. Observe that $-\m_1\l_1>-\m_2\l_1$ if $\l_1>0$ and $-\m_1\l_1<-\m_2\l_1$ if $\l_1<0$. Hence, we distinguish several cases:
\begin{enumerate}
\item Assume that $\m^*\leq \m_1<\m_2<0$. In this case, $h_1({\m_2})$ and $h_1(\m_1)$ are negative, and then we compare the functions $f_{\m_2}$ and $f_{\m_1}$ for negative values. Indeed, observe that $f_{\m_2}(\l_1)>f_{\m_1}(\l_1)$ for $\l_1<0$, and then $h_1({\m_2})<h_1(\m_1)$.
\item Assume that $\m_1<\m^*<\m_2<0$: in this case $h_1({\m_2})<0<h_1(\m_1)$.
\item Assume that $\m_1<\m_2\leq \m^*<0$, then $f_{\m_2}(\l_1)<f_{\m_1}(\l_1)$ for $\l_1>0$,
and then $h_1({\m_2})<h_1(\m_1)$.
\end{enumerate}   
This shows that $\m\mapsto h_1(\m) $ is decreasing.

We prove now the continuity. Take $\m_n\in (-\infty, 0)\to \m_0<0$ and consider $\l_n:=h_1({\m_n})$.  Since $0=f_{\m_n}(\l_n)=F(\l_n,\m_n\l_n)$,  by Corollary \ref{cotaraiz} we conclude that
\begin{equation}
\label{cotacaso1}
\l_n<\L_1^+\quad\mbox{and}\quad \m_n\l_n<\L_2^+.
\end{equation}
Hence, there exists $\overline{\l}_1\in (-\infty,+\infty)$ such that $\l_n\to \overline{\l}_1$. We have to show that 
$$
\overline{\l}_1=h_1(\m_0).
$$
Indeed, observe that
$$
0=f_{\m_n}(\l_n)=\L_1(-\l_n m_1,-\l_n \m_n m_2)\to \L_1(-\overline{\l}_1 m_1,-\overline{\l}_1 \m_0 m_2)=f_{\m_0}(\overline{\l}_1),
$$
that is, $f_{\m_0}(\overline{\l}_1)=0$.  We separate now two cases:
\begin{enumerate}
\item $\m_0\neq \m^*$: In this case, we assert. that $\overline{\l}_1\neq 0$. Indeed, assume that $\overline{\l}_1= 0$, that is, $h_1({\m_n})\to 0$. If, for instance,   $\m_0> \m^*$, then there exists $\rho_1({\m_n})\in (h_1({\m_n}),0)$ such that $f_{\m_n}'(\rho_1({\m_n}))=0$. Since $h_1({\m_n})\to 0$, then 
 $\rho_1({\m_n})\to 0$, and  as consequence, $f_{\m_0}'(0)=0$, a contradiction. This shows that $\overline{\l}_1\neq 0$. Then, since $h_1(\m_0)$ is the unique nonzero root of $f_{\m_0}(\l_1)=0$, we have that $\overline{\l}_1=h_1(\m_0)$.
\item $\m_0= \m^*$: in this case $f_{\m^*}(\overline{\l}_1)=0$ implies that $\overline{\l}_1=0=h_1({\m^*})$.
\end{enumerate}
This concludes that  $\overline{\l}_1=h_1({\m_0})$, and hence the continuity.

We claim that
\begin{equation}
\label{claim1}
h_1(\m_n)\to -\infty\quad\mbox{as $\m_n\to 0$.}
\end{equation} 
Assume that $|h_1(\m_n)|\leq C$. Then, we can assume that, at least for a subsequence, $h_1(\m_n)\to h_1^*<0$ and hence
$$
0=\L_1(-h_1(\m_n)m_1,-\m_nh_1(\m_n)m_2)\to \L_1(-h_1^*m_1,0)=0,
$$ 
a contradiction because $ \L_1(-h_1^*m_1,0)>0$. This proves (\ref{claim1}).

By  (\ref{cotacaso1}), if $\m\to -\infty$ we can assume that $h_1(\m)\to h^*\leq \L_1^+$ and $h^*>0$. Then, $\m h_1(\m)\to -\infty$. Since
$$
0=f_{\m}(h_1(\m))=F(h_1(\m),h_1(\m)\m)
$$
and by Proposition  \ref{convergence}
$$
0=F(h_1(\m),h_1(\m)\m)\to \min\{\s_1^{\O_1}(-\D-\l^*m_1;N+\g_1),\s_1^{M_2^0}(-\D;D)\},
$$
it follows that
$$
0= \min\{\s_1^{\O_1}(-\D-\l^*m_1,N+\gamma_1),\s_1^{M_2^0}(-\D,D)\}.
$$
Since $\s_1^{M_2^0}(-\D,D)>0$, we conclude that $\s_1^{\O_1}(-\D-h^*m_1,N+\gamma_1)=0$, 
that means that $h^*=\L_1^+$, that is
$$
\lim_{\mu\to -\infty}h_1(\m)=\L_1^+.
$$
 \end{enumerate}
This concludes the proof. 
 \end{proof}

\begin{figure}[h]
\centering
\includegraphics[scale=.25]{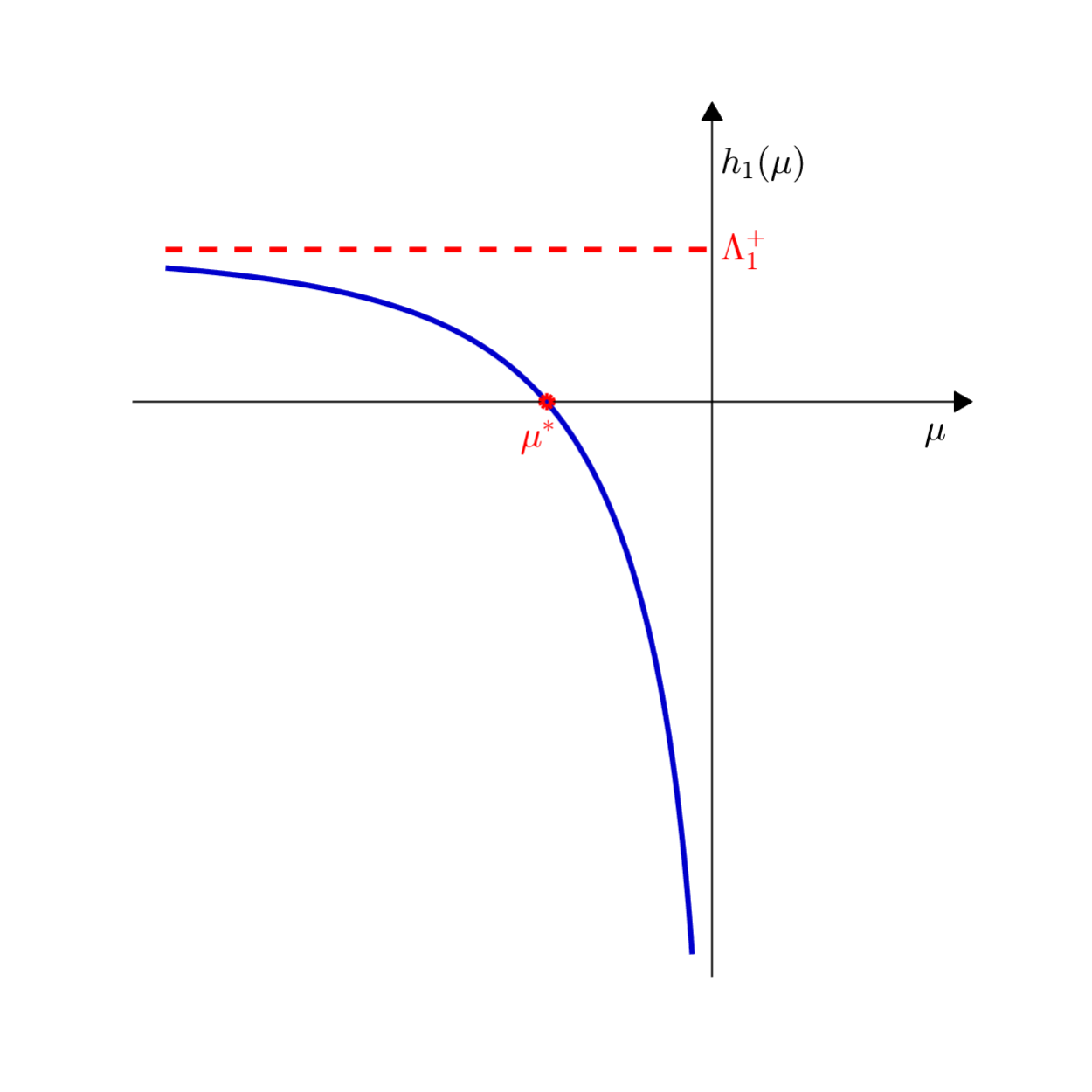}
\includegraphics[scale=.25]{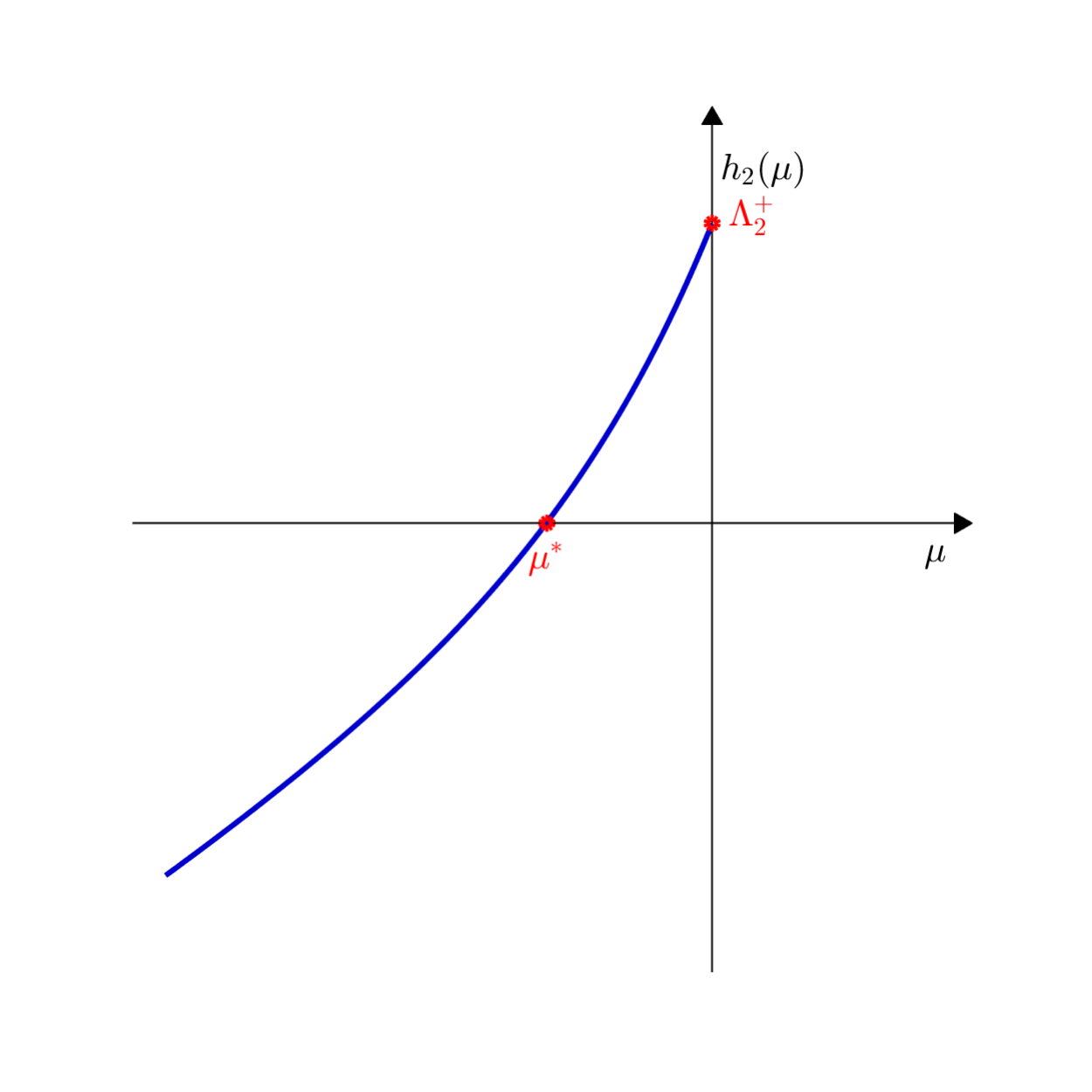}
\caption{Case $m_i\gneq 0$ in $\O_i$ and verifying (\ref{00intro}) for $i=1,2$. We have represented the functions $\m\mapsto h_1(\m)$ (left) and $\m\mapsto h_2(\m)$ (right).}
\label{figuram1}
\end{figure}

Once we have studied the map $\m\mapsto h_1(\m)$, we need to analyze the map 
$$\m\in (-\infty,0)\mapsto h_2(\m):=\m h_1(\m).
$$  
 \begin{proposition}
 \label{lambda2}
 Assume that $m_i\gneq 0$ in $\O_i$ and $M_i^0$ verify (\ref{00intro}) for $i=1,2$.   The map $\m\in (-\infty,0)\mapsto h_2(\m):=\m h_1(\m)$ is continuous, increasing,
 $$
 h_2(\m)
 \left\{
 \begin{array}{ll}
 >0 & \mbox{if $\m>\m^*$ },\\
 =0 & \mbox{if $\m=\m^*$,}\\
 <0 & \mbox{if $\m<\m^*$,}
 \end{array}  
 \right.
 $$
 \begin{equation}
 \label{marka}
 \lim_{\m\to -\infty}h_2(\m)=-\infty,
  \end{equation}
and 
  \begin{equation}
 \label{markb}
 \lim_{\m\to 0}h_2(\m)=\L_2^+.
 \end{equation}
 \end{proposition}
 \begin{proof}
 To start with, the continuity and the sign of the map $h_2(\m)$ follow directly from Proposition~\ref{casegeneraligualsigno}. Moreover, it is clear that
 $$
 \lim_{\m\to -\infty}h_2(\m)=\lim_{\m\to -\infty}\m h_1(\m)=-\infty.
 $$
 In order to prove (\ref{markb}) we can argue exactly as the proof of Proposition~\ref{casegeneraligualsigno}. 
 
 Finally, using that $F$ is increasing, we prove that the map $\m\mapsto h_2(\m)$ is increasing. Take $\m_1<\m_2$ and assume that $h_2(\m_1)\geq h_2(\m_2)$. Since $h_1(\m_1)>h_1(\m_2)$, then 
 $$
 0=F(h_1(\m_1),h_2(\m_1))<F(h_1(\m_2),h_2(\m_1))\leq F(h_1(\m_2),h_2(\m_2))=0,
 $$
 a contradiction.
 \end{proof}
  In Figure \ref{figuram1} we have represented the functions $\m\mapsto h_1(\m), h_2(\m)$. Now, we proceed to the proof of Theorem \ref{teoprinci1}.
  
 \begin{proof}[Proof of Theorem \ref{teoprinci1}:] (see Figure \ref{figura:figure1}) 
 
 \begin{enumerate}
\item  Observe that by Corollary \ref{cotaraiz}, we obtain
 $$
F(\l_1,\l_2)<0\quad\mbox{if $\l_1\geq \L_1^+$ or $\l_2\geq \L_2^+$.}
$$
 
\item Take $\l_1<\L_1^+$. Then, by Proposition \ref{casegeneraligualsigno} there exists a unique $\m=\m(\l_1)<0$ such that $\l_1=h_1(\m)$.  Take $\l_2=h_2(\m)=\m h_1(\m)$, then
$$
F(\l_1,\l_2)=F(h_1(\m),h_2(\m))=F(h_1(\m),\m h_1(\m))=f_\mu(h_1(\m))=0.
$$
We define the function
$$
{\cal H}(\l_1):=h_2(h_1^{-1}(\l_1)),\qquad \l_1<\L_1^+.
$$
It is clear that ${\cal H}$ is well-defined (observe that $h_1^{-1}$ exists due to that $h_1$ is a decreasing function), is continuous and 
$$
F(\l_1,{\cal H}(\l_1))=0.
$$
Moreover, since fixed $\l_1$, the map $\l_2\mapsto F(\l_1,\l_2)$ is concave, it follows that $F(\l_1,\l_2)>0$ for $\l_2<{\cal H}(\l_1)$ and  $F(\l_1,\l_2)<0$ for $\l_2>{\cal H}(\l_1)$.

Furthermore, by Propositions \ref{casegeneraligualsigno}  and \ref{lambda2},
$$
\lim_{\l_1\uparrow \L_1^+}{\cal H}(\l_1)=\lim_{\l_1\uparrow \L_1^+}h_2(h_1^{-1}(\l_1))= \lim_{\m\to -\infty}h_2(\m)=-\infty,
$$
and, 
$$
 \lim_{\l_1\uparrow  -\infty}{\cal H}(\l_1)= \lim_{\l_1\uparrow  -\infty}h_2(h_1^{-1}(\l_1))=\lim_{\m\to 0}h_2(\m)=\L_2^+.
$$
Finally, we prove that $\l_1\mapsto {\cal H}(\l_1)$ is decreasing. Take $\l_1<\overline{\l}_1<\L_1^+$ and consider ${\cal H}(\l_1)$ and ${\cal H}(\overline{\l}_1)$. We are going to show that ${\cal H}(\l_1)>{\cal H}(\overline{\l}_1)$. Assume that ${\cal H}(\l_1)\leq {\cal H}(\overline{\l}_1)$, then
$$
0=F(\l_1,{\cal H}(\l_1))>F(\overline{\l}_1,{\cal H}(\l_1))\geq F(\overline{\l}_1,{\cal H}(\overline{\l}_1))=0,
$$
a contradiction.
\end{enumerate}
This concludes the proof.
\end{proof}
\subsection{Case $m_1\gneq 0$ in $\O_1$ and $m_2$ changes sign in $\O_2$.}
In this case, the results depend on the sign of $\int_{\O_2}m_2$. We detail the case
$$
\int_{\O_2}m_2<0,
$$
similarly the other cases can be studied (see Remark \ref{casoparticular1}). Observe that in this case
$$
\m^*=-\frac{\g_2\displaystyle\int_{\O_1}m_1}{\g_1\displaystyle\int_{\O_2}m_2}>0.
$$
\begin{proposition}
\label{posics}
Assume that $m_1\gneq 0$ in $\O_1$ and $M_1^0$ verifies (\ref{00intro}), $m_2$ changes sign in $\O_2$ and  $
\int_{\O_2}m_2<0.$ Then, for each $\m\neq 0$ there exists a unique $h_1(\m)\in \R$ such that $f_\m(h_1(\m))=0$. Moreover, 
 $$
 h_1(\m)
 \left\{
 \begin{array}{ll}
 >0 & \mbox{if $\m>\m^*$ },\\
 =0 & \mbox{if $\m=\m^*$,}\\
 <0 & \mbox{if $\m<\m^*$.}
 \end{array}  
 \right.
 $$
Furthermore, the map $\m\in \R\setminus\{0\}\mapsto h_1(\m)\in\R$ is continuous, and
$$
\lim_{\m\to \pm\infty}h_1(\m)=0,\quad\lim_{\m\to 0}h_1(\m)=-\infty.
$$
As consequence, there exists $\m_{\max}>\m^*$ such that
$$
\max_{\m\neq 0}h_1(\m)=h_1(\m_{\max}):=\l_1^{max}.
$$
Finally, the map $\m\in \R\setminus\{0\}\mapsto h_1(\m)\in\R$ is increasing in $(0,\m_{\max})$ and decreasing in $(-\infty,0)$ and $(\m_{\max},\infty)$.
\end{proposition}
\begin{proof}
The proof of this result is rather similar to the proof of Proposition \ref{casegeneraligualsigno}, hence we sketch the proof. 

Since $f_{\m_n}(\l_1(\m_n))=F(h_1(\m_n),\m_n h_1(\m_n))=0$, by Corollary \ref{cotaraiz} we get
$$
\L_2^-<h_1(\m_n)\m_n<\L_2^+,
$$
whence we conclude that $h_1(\m_n)\to 0$ as $\m_n\to \pm\infty$.

Before proving the monotony, we claim that for any $c\in\R$ there exist at most two values of $\m$  such that
$$
h_1(\m)=c.
$$  
We argue by contradiction. Assume that for $\m_1<\m_2<\m_3$ we get $h_1(\m_i)=c$ for $i=1,2,3$. Taking $\l_2^i= c\m_i$ we obtain
$$
0=F(c,\l_2^i),\quad \l_2^1<\l_2^2<\l_2^3,
$$
a contradiction because, fixed $c$, the map $\l_2\mapsto F(c,\l_2)$ is concave. 

Now, for instance, we show that $h_1(\m)$ is decreasing in $(-\infty,0)$. Take $\m_1<\m_2<0$ and assume that $h_1(\m_1)\leq h_1(\m_2)$. Since $h_1(\m)\to 0$ as $\m\to -\infty$ and $h_1(\m)\to -\infty$ as $\m\to 0$. Then, there exists $c<0$ such that $h_1(\m)=c$ possesses at least three solutions. This is a contradiction and proves that $h_1(\m)$ is decreasing in $(-\infty,0)$.  With a similar argument, it can be proved that $h_1(\mu)$ is decreasing in  $(\m_{\max},\infty)$ and increasing in  $(0,\m_{\max})$ .
 \end{proof}

Again, we can deduce properties of the map $h_2(\m)=\m h_1(\m)$. 
\begin{proposition}
\label{posicsl2}
Assume that $m_1\gneq 0$ in $\O_1$ and $M_1^0$ verifies (\ref{00intro}), $m_2$ changes sign in $\O_2$ and  $
\int_{\O_2}m_2<0.$  Then $h_2(\m)=\m h_1(\m)$ is continuous in $\m\neq 0$, increasing and  verifies
 $$
 h_2(\m)
 \left\{
 \begin{array}{ll}
 >0 & \mbox{if $\m>\m^*$ or $\m<0$},\\
 =0 & \mbox{if $\m=\m^*$,}\\
 <0 & \mbox{if $\m\in (0,\m^*)$,}
 \end{array}  
 \right.
 $$
and 
$$
\lim_{\m\to 0^\pm}h_2(\m)=\L_2^\mp, \quad \lim_{\m\to\pm\infty}h_2(\m)=\l_2^*,
$$
for some $\l_2^*\in(0,\L_2^+)$. 
\end{proposition}
\begin{proof}
Assume that $\m_n\to 0^+$, then since $h_2(\m_n)$ is bounded, at least for a subsequence, $h_2(\m_n)\to \overline{\l_2}<0$. Observe that
since $h_1(\m_n)\to -\infty$, by Proposition \ref{convergence} 
$$
0=F(h_1(\m_n),h_2(\m_n) )\to \s_1^{\O_2}(-\D-\overline{\l_2}m_2;N+\g_2,N)=0,
$$ 
and then $\overline{\l_2}=\L_2^-$. Analogously for $\m_n\to 0^-$.

On the other hand, assume that $\m_n\to +\infty$ and $h_2(\m_n)\to \overline{\l_2}<0$. In this case, $h_1(\m_n)\to 0$, and then
$$
0=F(-h_1(\m_n)m_1,-h_2(\m_n)m_2)\to F(0,-\overline{\l_2}m_2),
$$ 
whence $\overline{\l_2}=\l_2^*.$

Observe that 
$$
0=F(0,-\l_2^*m_2)<\s_1^{\O_2}(-\D-\l_2^*m_2;N+\g_2,N)
$$
and so $\l_2^*<\L_2^+$.

Finally, with an argument similar to the one used in Proposition  \ref{posics} we can conclude that the equation $h_2(\m)=c$ possesses at most a unique solution. Hence, the monotony of $h_2(\m)$ follows. This completes the proof. 
\end{proof}

\begin{figure}[h]
\centering
\includegraphics[scale=.35]{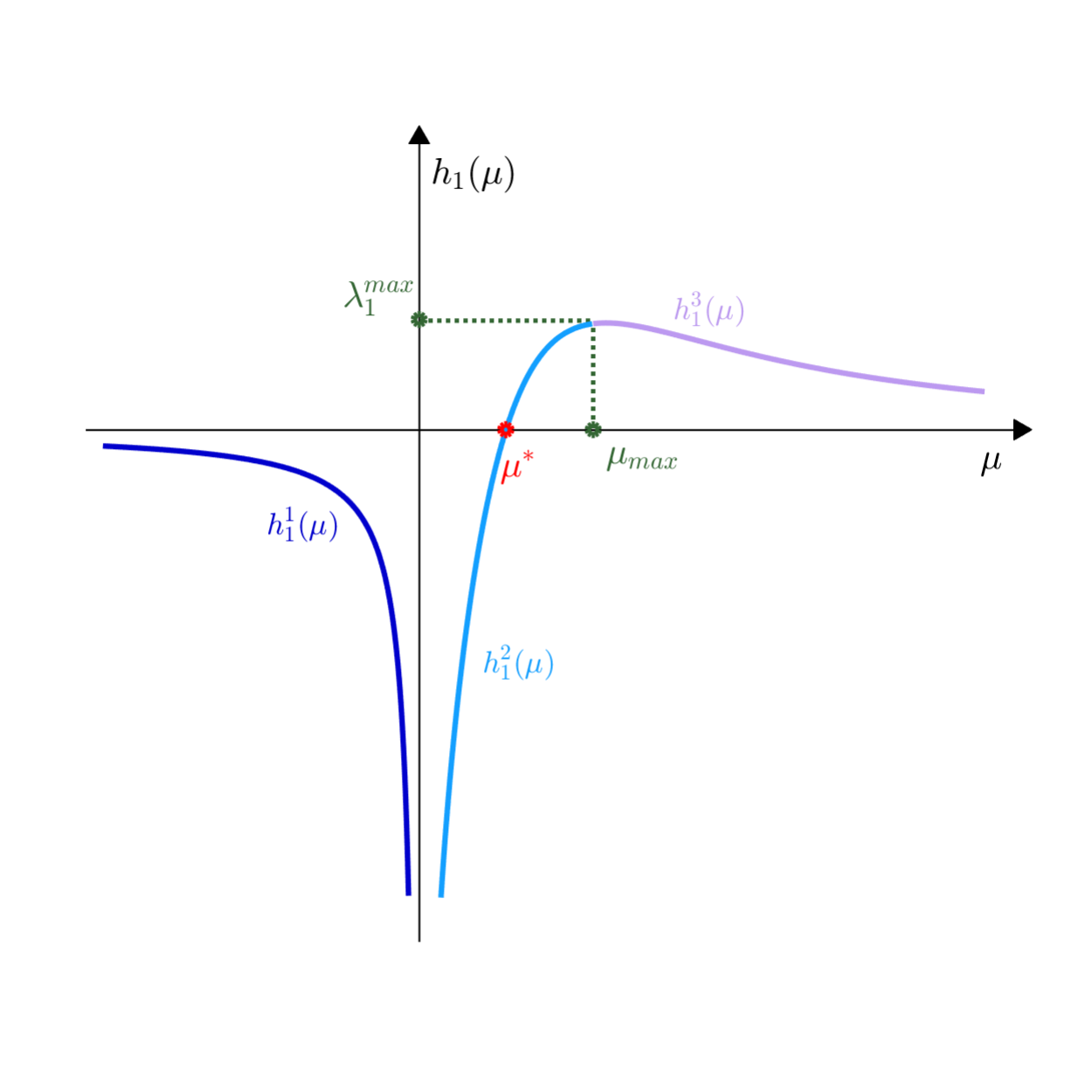}
\includegraphics[scale=.35]{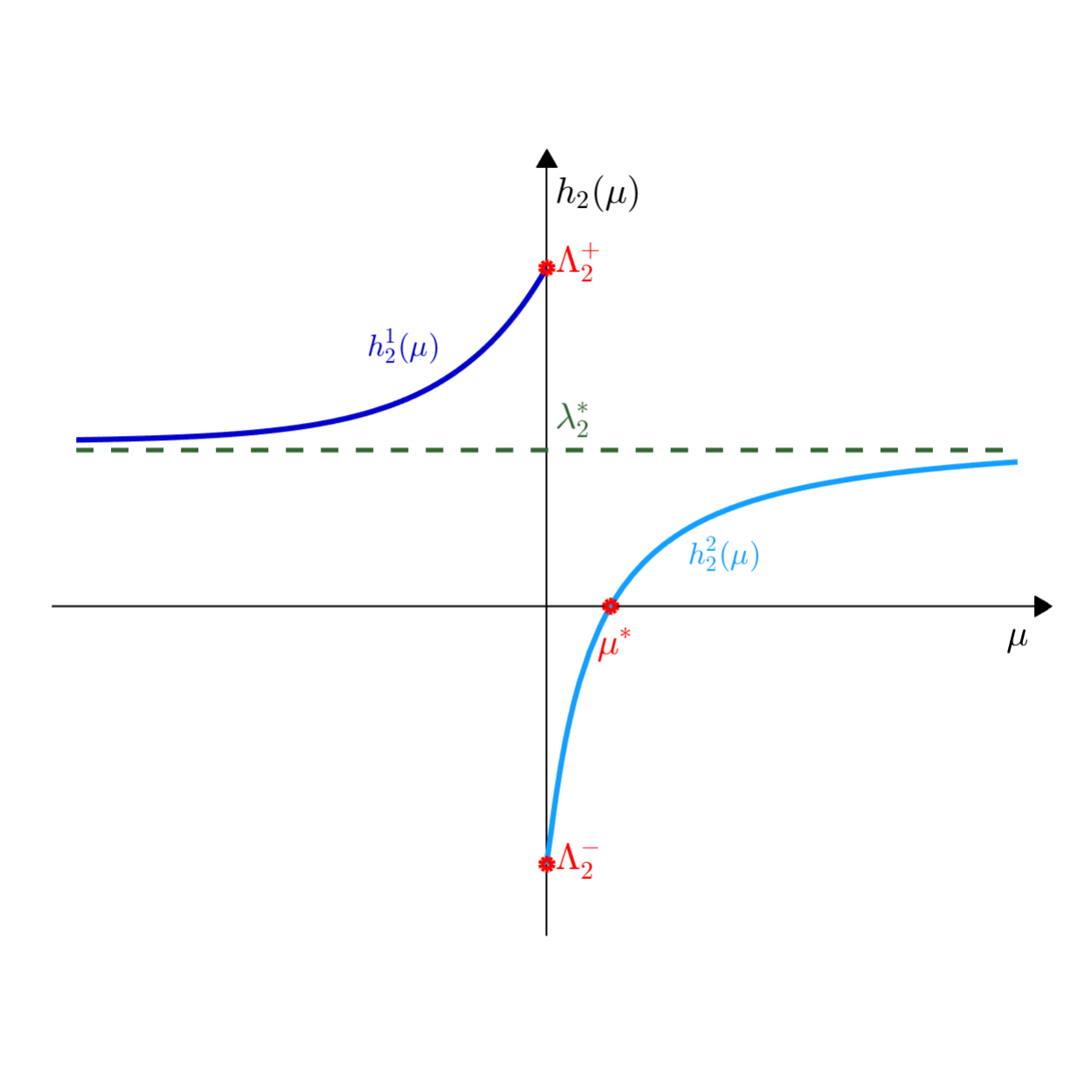}
\caption{Representation of the maps $h_1(\m)$ (left)  and $h_2(\m)$ (right)  in the case $m_1$ non-negative, non-trivial and verifies (\ref{00intro}), $m_2$ changing sign and $\int_{\O_2}m_2<0$.}
\label{figuram2}
\end{figure}

In Figure \ref{figuram2}, one may see a representation of the maps $\m\mapsto h_1(\m),h_2(\m)$.

\begin{figure}[h]
\centering
\includegraphics[scale=.30]{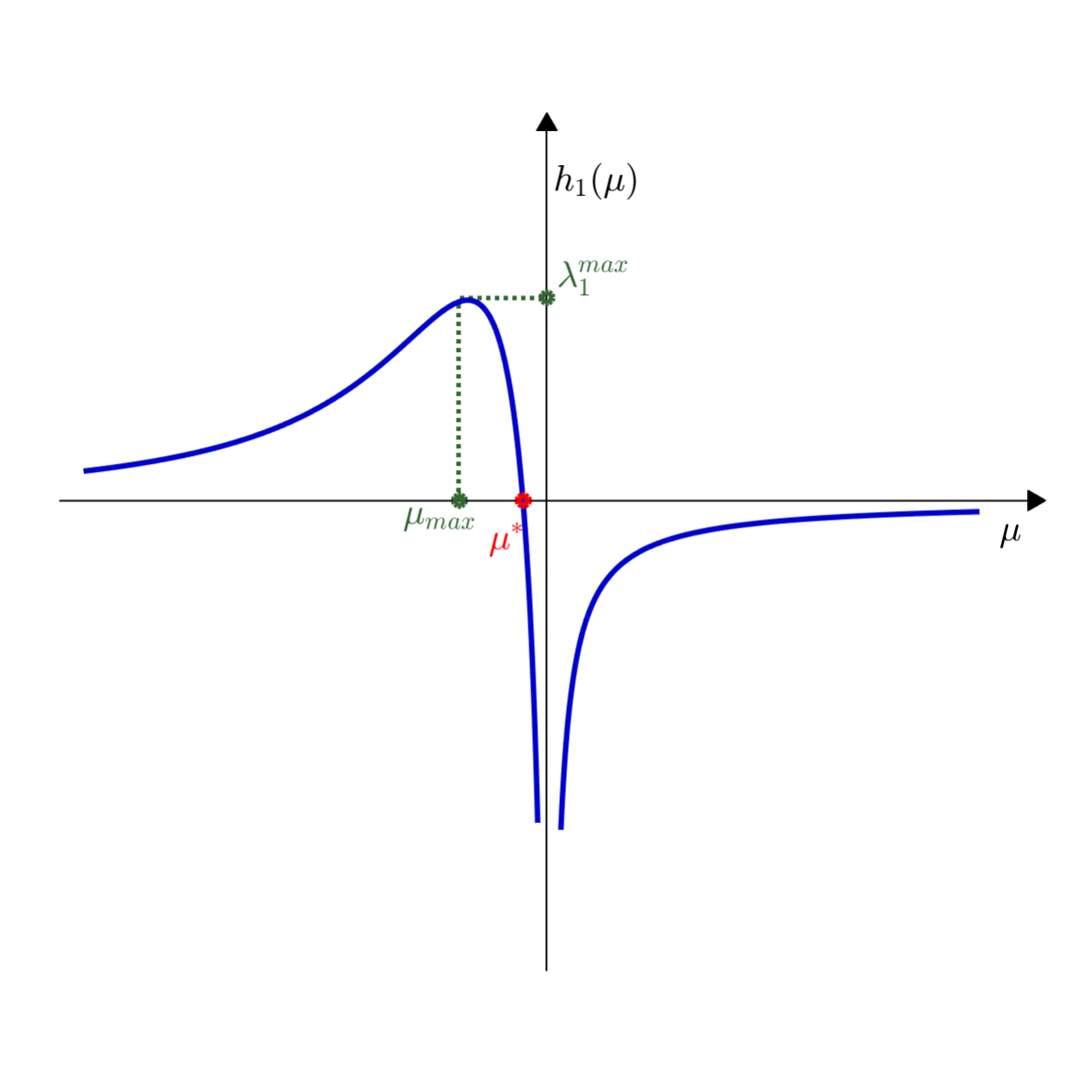}
\includegraphics[scale=.30]{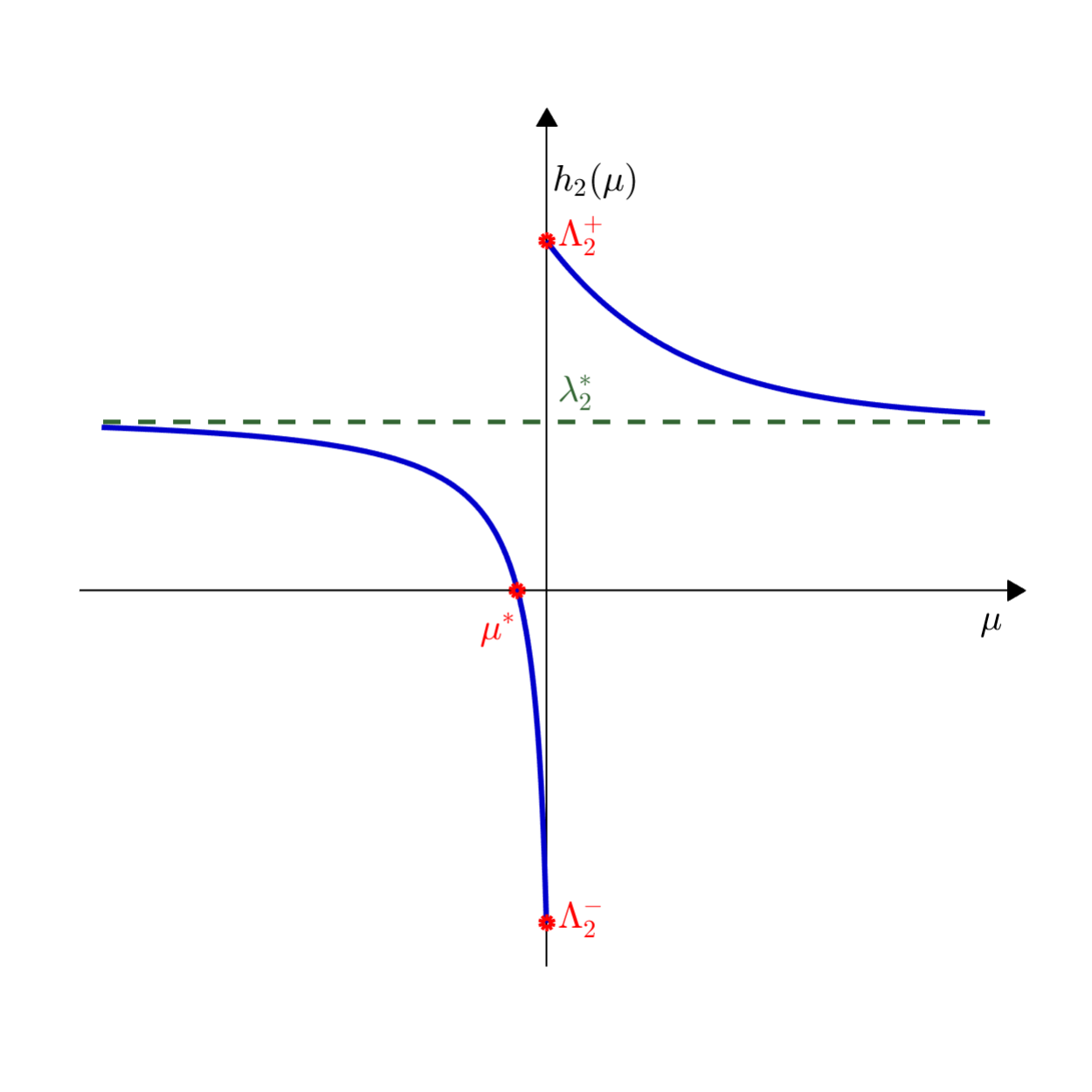}
\caption{Representations of the maps $h_1(\m)$ (left)  and $h_2(\m)$ (right)  in the case $m_1$ non-negative, non-trivial and verifies (\ref{00intro}), $m_2$ changing sign and $\int_{\O_2}m_2>0$.}
\label{figuram3}
\end{figure}

\begin{figure}[h]
\centering
\includegraphics[scale=.30]{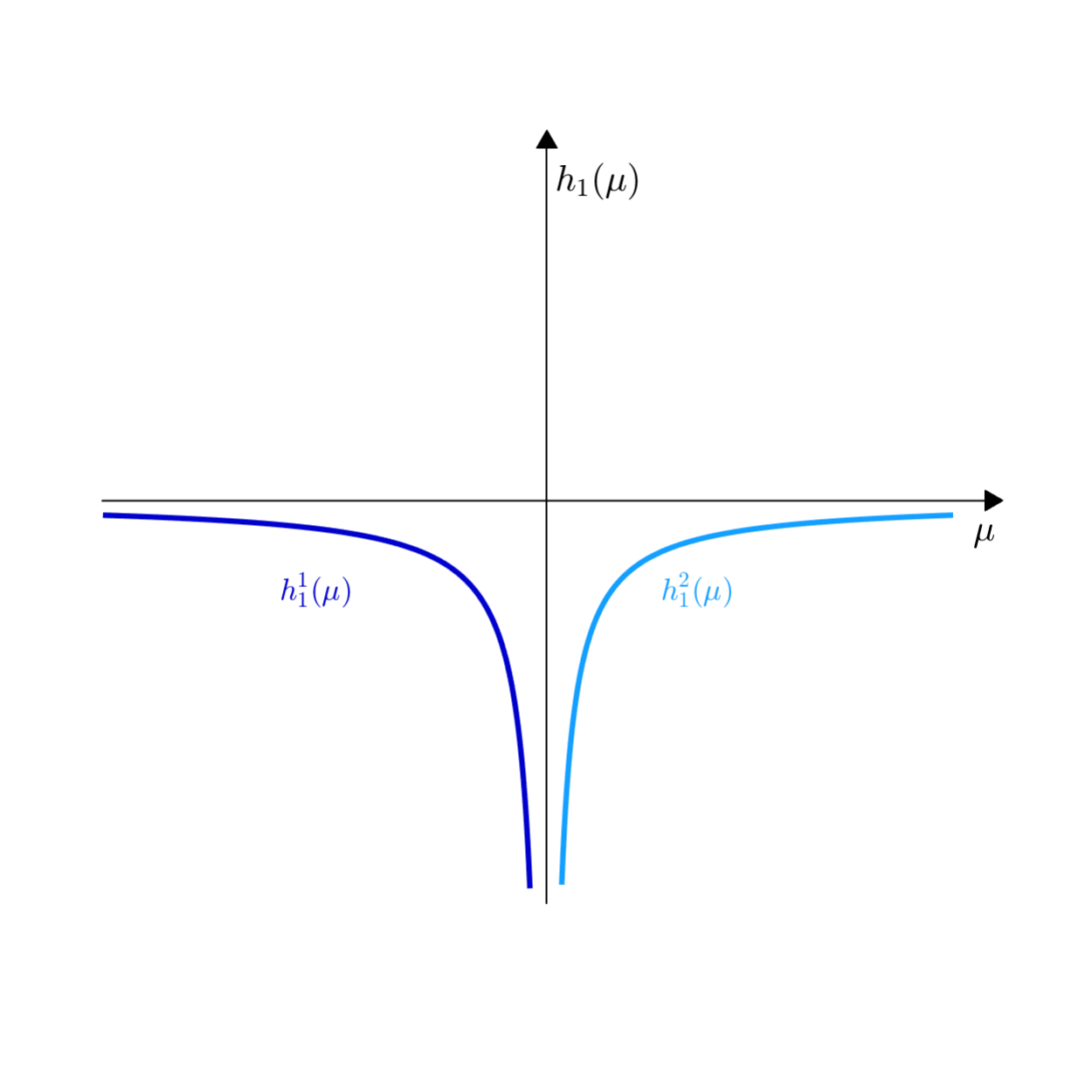}
\includegraphics[scale=.30]{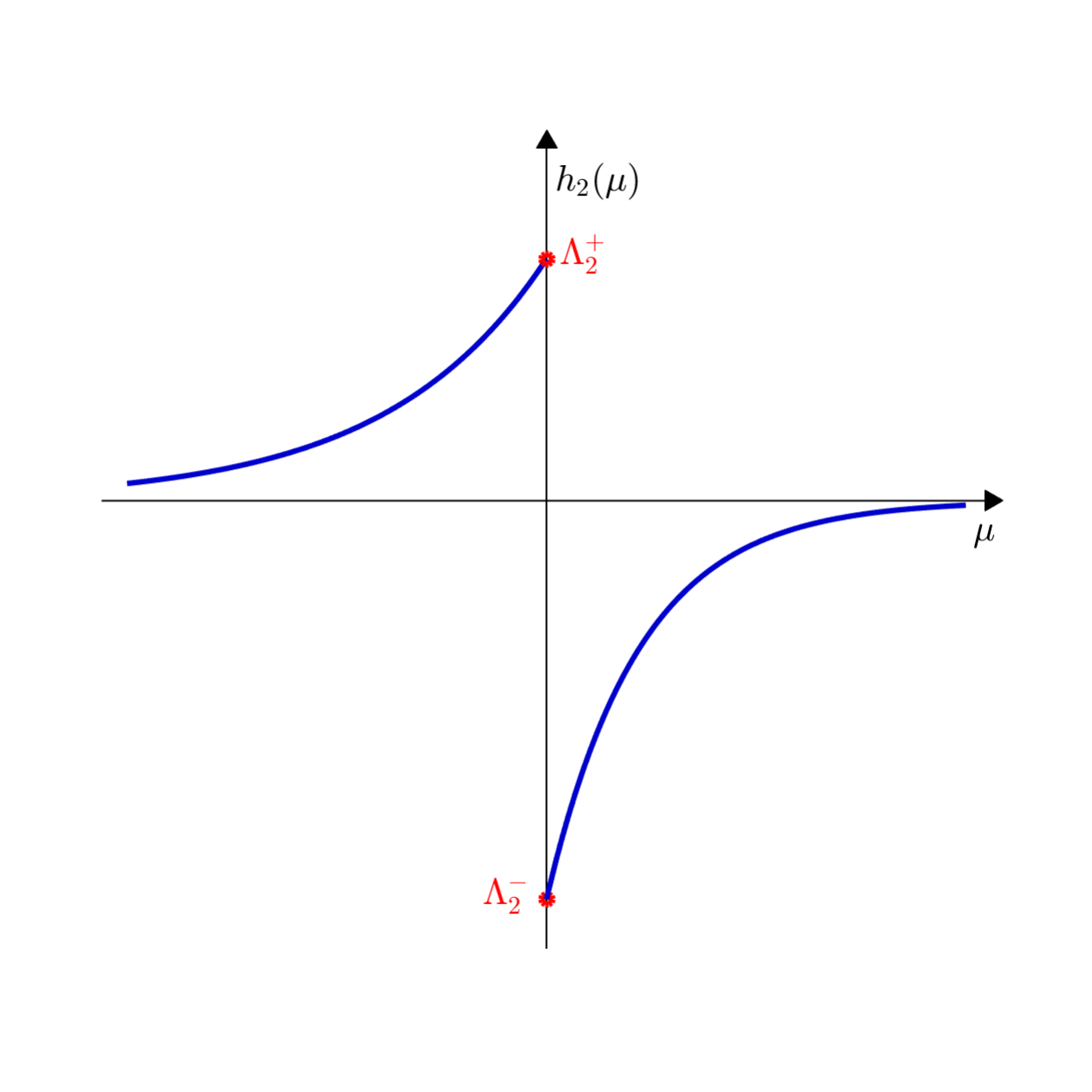}
\caption{We have represented the maps $h_1(\m)$ (left)  and $h_2(\m)$ (right)  in the case $m_1$ non-negative, non-trivial and verifies (\ref{00intro}), $m_2$ changing sign and $\int_{\O_2}m_2=0$.}
\label{figuram4}
\end{figure}

\begin{remark}
\label{casoparticular1}
\begin{enumerate}
\item  In the case 
$$
\int_{\O_2}m_2>0
$$
we can obtain a similar result switching $\m$ by $-\m$ (see Figure \ref{figuram3}).
\item When 
$$
\int_{\O_2}m_2=0
$$
observe that $f'_\mu(0)<0$ for all $\m$ (see (\ref{derivada})), and then $h_1(\m)<0$ for all $\m$. The global behavior of $h_1(\m)$ at $\m=0$ and $\m\to \pm\infty$ is  similar to Proposition \ref{posics} (see Figure  \ref{figuram4}).
\end{enumerate}
\end{remark}

\begin{proof}[Proof of Theorem \ref{teoprinci2}:] (See Figure \ref{fig:sub1}).  Assume that $\int_{\O_2}m_2<0$ (see Figure \ref{figuram2}).
We introduce the following notation:
$$
h_1(\m):=
\left\{
\begin{array}{ll}
h_1^1(\m) & \mbox{if $\m<0$,}\\
h_1^2(\m) & \mbox{if $\m\in (0,\m_{\max}$],}\\
h_1^3(\m) & \mbox{if $\m>\m_{\max}$.}
\end{array}
\right.
$$
%and 
%$$
%h_2(\m):=
%\left\{
%\begin{array}{ll}
%h_2^1(\m) & \mbox{if $\m<0$,}\\
%h_2^2(\m) & \mbox{if $\m>0$.}
%\end{array}
%\right.
%$$

\begin{enumerate}
\item If $\l_1>\l_1^{\max}$, then there does not exist $\m\in\R$ such that $\l_1=h_1(\m)$. Hence,  $F(\l_1,\l_2)\neq0$ for all $\l_2\in\R$, in fact,  $F(\l_1,\l_2)<0$ for all $\l_2\in\R$. Indeed, if for some $\overline{\l_2}$ we have $F(\l_1,\overline{\l_2})>0$, then there exists at least $\l_2^0$ such that $F(\l_1,\l_2^0)=0$. Then, for some $\m_0$ we have $\l_1=h_1(\m_0)$, a contradiction.

\item  If $\l_1=\l_1^{\max}$, there exists a unique $\m_{max}>\m^*$ such that $\l_1^{\max}=h_1(\m_{\max})$, and then $\l_2^{\max}=h_2(\m_{\max})>0$ and $F(\l_1^{\max},\l_2^{\max})=0$.

\item We fix $\l_1\in (0,\l_1^{\max})$. Then, (see Figure \ref{figuram2}) there exist $\m_2,\m_3$ with $\m^*<\m_2<\m_{max}<\m_3$ such that $h_1^i(\m_i)=\l_1$ $i=2,3$. To these values correspond two different values of $h_2(\m_i)$. Moreover, as $\l_1\to 0$, then $\m_2\to \m^*$ and $\m_3\to +\infty$, and this case $h_2(\m_2)=h_2^2(\m_2)\to h^2_2(\m^*)=0$ and $h_2(\m_3)=h_2^2(\m_3)\to \l_2^*$.

\item On the other hand, when $\l_1\in (-\infty,0)$. There exist $\m_1<0<\m_2<\m^*$ such that $\l_1=h_1^i(\m_i)$ $i=1,2,$ with $\m_1\to -\infty$ and $\m_2\to \m^*$ as $\l_1\to 0$. Then, $h_2(\m_1)=h_2^1(\m_1)\to \l_2^*$ and $h_2(\m_2)=h_2^2(\m_2)\to 0$.   

Observe that when $\l_1\to -\infty$ then $\m_1\to 0^-$ and $\m_2\to 0^+$, and hence $h_2(\m_1)\to \L_2^+$ and  $h_2(\m_2)\to \L_2^-$.
\end{enumerate}
With this construction, we can define
$$
{\cal H}^+(\l_1):=
\left\{
\begin{array}{ll}
h_2((h_1^3)^{-1}(\l_1)) & \mbox{if $\l_1\in (0,\l_1^{\max}],$}\\
h_2((h_1^1)^{-1}(\l_1)) & \mbox{if $\l_1\leq 0,$}\\
\end{array}
\right.
$$
and
$$
{\cal H}^-(\l_1):=h_2((h_1^2)^{-1}(\l_1))\quad\mbox{if $\l_1\leq\l_1^{\max}.$}
$$
Observe that thanks to the monotony of the maps $h_1^i$ for $i=1,2,3$, ${\cal H}^+$ and ${\cal H}^-$ are well defined. Moreover,
$$
\lim_{\l_1\to 0^+}{\cal H}^+(\l_1)=\lim_{\l_1\to 0^+}h_2((h_1^3)^{-1}(\l_1)) =\lim_{\m\to +\infty}h_2(\m)=\l_2^*,
$$ 
and 
$$
\lim_{\l_1\to 0^-}{\cal H}^+(\l_1)=\lim_{\l_1\to 0^-}h_2((h_1^1)^{-1}(\l_1)) =\lim_{\m\to -\infty}h_2(\m)=\l_2^*.
$$ 
As a consequence, ${\cal H}^+$ is continuous. 

On the other hand,
$$
\lim_{\l_1\to\l_1^{\max}}{\cal H}^+(\l_1)=\lim_{\m\to\m_{\max}}h_2(\m)=\overline{\l}_2,
$$
$$
\lim_{\l_1\to\l_1^{\max}}{\cal H}^-(\l_1)=\lim_{\m\to\m_{\max}}h_2(\m)=\overline{\l}_2.
$$
Finally,
$$
\lim_{\l_1\to-\infty}{\cal H}^+(\l_1)=\lim_{\m\to 0^+}h_2(\m)=\L_2^+,
$$
and
$$
\lim_{\l_1\to-\infty}{\cal H}^-(\l_1)=\lim_{\m\to 0^-}h_2(\m)=\L_2^-.
$$

We show that ${\cal H}^+(\l_1)$ is decreasing. Take $\l_1^1<\l_1^2$. 
\begin{enumerate}
\item When $\l_1^1<\l_1^2<0$: then $(h_1^1)^{-1}(\l_1^2)<(h_1^1)^{-1}(\l_1^1)<0$ and so $h_2((h_1^1)^{-1}(\l_1^2))<h_2((h_1^1)^{-1}(\l_1^1)).$ This concludes that
$$
{\cal H}^+(\l_1^1)>{\cal H}^+(\l_1^2).
$$
\item Assume now that $\l_1^1<0<\l_1^2$:  in this case $(h_1^1)^{-1}(\l_1^1)<0< (h_1^3)^{-1}(\l_1^2)$ and then $h_2((h_1^1)^{-1}(\l_1^1))>0>h_2((h_1^3)^{-1}(\l_1^2) )$, that is ${\cal H}^+(\l_1^1)>{\cal H}^+(\l_1^2).$
\item Finally when $0<\l_1^1<\l_1^2$: in this case $0< (h_1^3)^{-1}(\l_1^2)<(h_1^3)^{-1}(\l_1^1)$.  Again,  ${\cal H}^+(\l_1^1)>{\cal H}^+(\l_1^2)$.
\end{enumerate}
We can argue in the same manner for ${\cal H}^-$. This completes the proof.
\end{proof}
\begin{remark}
\label{casem0}
Case $\int_{\O_2}m_2>0$ can be handled in an analogous way, but the case $\int_{\O_2}m_2=0$ deserves a comment. In this case, ${\cal H}^+$ and ${\cal H}^-$ should be defined as follows:
$$
{\cal H}^+(\l_1):=
\left\{
\begin{array}{ll}
h_2((h_1^1)^{-1}(\l_1)) & \mbox{if $\l_1<0,$}\\
0 & \mbox{if $\l_1=0$,}\\
\end{array}
\right.
$$
and
$$
{\cal H}^-(\l_1):=
\left\{
\begin{array}{ll}
h_2((h_1^2)^{-1}(\l_1)) & \mbox{if $\l_1<0,$}\\
0 & \mbox{if $\l_1=0$.}\\
\end{array}
\right.
$$

\end{remark}

%ESTO ES OTRA PRUEBA
%
%-----------------------------------------------------------------------------------------------
%
%In this case, 
%
%
%Again, we divide our study:
%\begin{enumerate}
%\item Assume $m>m^*$: in this case $f_m'(0)>0$ and so there exists $\l_1^m>0$ such that $f_m(\l_1^m)=0$. In this case, $\l_2^m=m\l_1^m>0$. Moreover, it can be shown that $\l_1^m\to 0$ and $\l_2^m\to \l_2^+$ as $m\to +\infty$. 
%
%Hence, we can define
%$$
%\inf_{m\geq m^*}\l_1m=\l_1^{m^{**}}>0.
%$$
%\item Assume $m=m^*$: no solution.
%\item Assume $0<m<m^*$: in this case, $\l_1^m<0$ and $\l_2^m<0$. Moreover,
%$$
%\l_1^m\to -\infty\quad\mbox{and}\quad \l_2^m\to\l_2^-\quad \mbox{as $m\to 0$}.  
%$$ 
%\item Assume $m<0$: In this case $\l_1^m<0$ and $\l_2^m>0$. Hay que ver lo que pasa con $m\to -\infty$. 
%\end{enumerate}
%-----------------------------------------------------------------------------------------------
\subsection{Case $m_i$ changes sign, $i=1,2$.}   
Consider in this case 
$$
\int_{\O_1}m_1<0,\quad \int_{\O_2}m_2<0,
$$
and then
$$
\m^*<0.
$$
\begin{proposition}
\label{l1cscs}
Assume that $m_i$ changes sign for $i=1,2$ and $
\int_{\O_1}m_1<0,\; \int_{\O_2}m_2<0$. Then, for each $\m\in\R$ there exists a unique value $\l_1=h_1(\m)$ such that $f_\m(\l_1)=0$.
Moreover,  the map $\m\in\R\mapsto h_1(\m)$ is continuous, 
$$
 h_1(\m)
 \left\{
 \begin{array}{ll}
 >0 & \mbox{if $\m>\m^*$},\\
 =0 & \mbox{if $\m=\m^*$,}\\
 <0 & \mbox{if $\m<\m^*$,}
 \end{array}  
 \right.
 $$
 and
$$
\lim_{\m\to \pm\infty}h_1(\m)=0.
$$
As consequence, there exist $\m_{\min}<\m^*<\m_{\max}$ such that
$$
h_1(\m_{\min})=\min_{\m\in\R} h_1(\m)=\l_1^{\min}<0,\quad h_1(\m_{max})=\max_{\m\in\R} h_1(\m)=\l_1^{\max}>0.
$$
Finally, the map $\m\mapsto h_1(\m)$ is decreasing in $(-\infty, \m_{\min})$ and $(\m_{\max},\infty)$  and increasing in $(\m_{\min},\m_{\max})$.
\end{proposition}

For $h_2(\m)$, we can deduce the following

\begin{figure}[h]
\centering
\includegraphics[scale=.35]{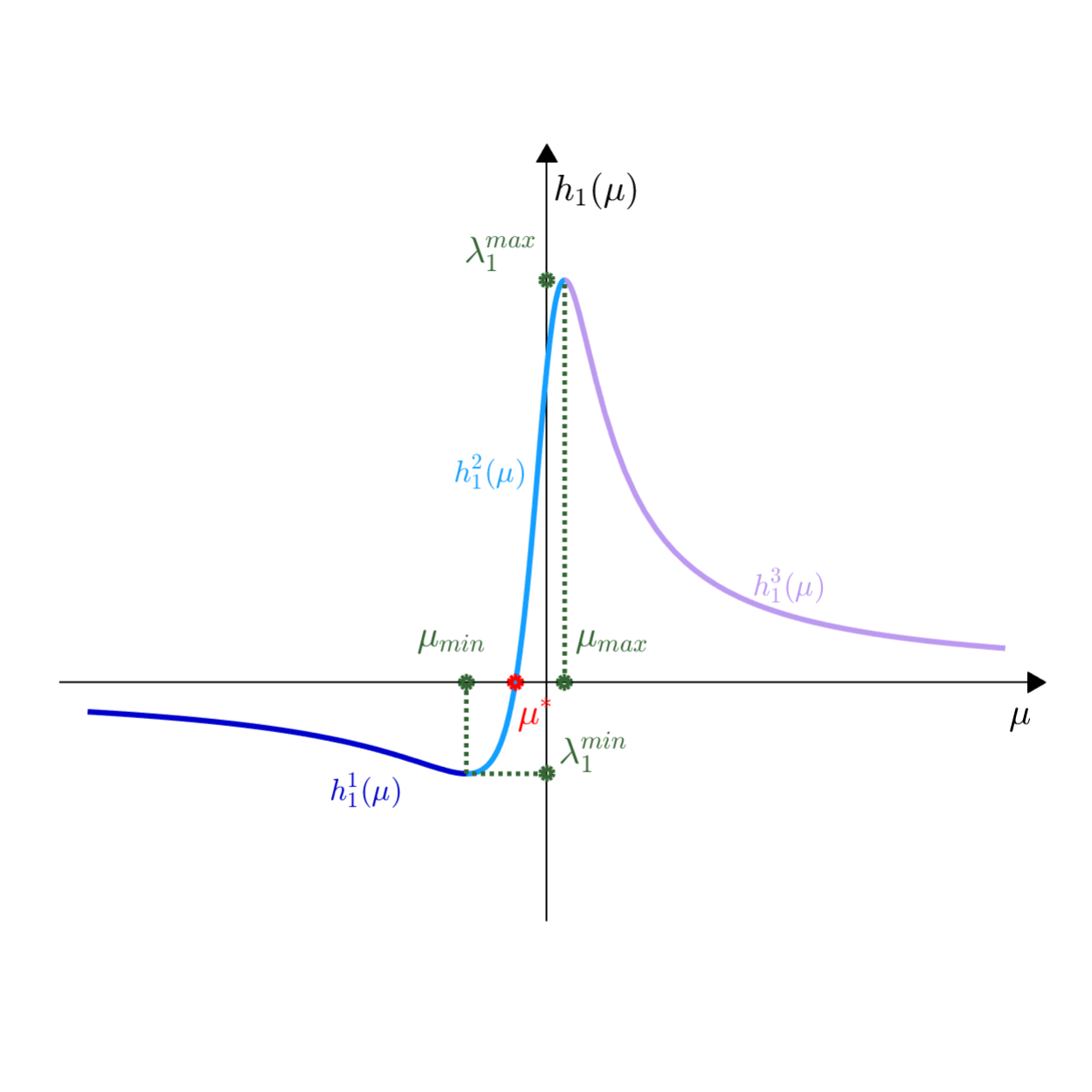}
\includegraphics[scale=.35]{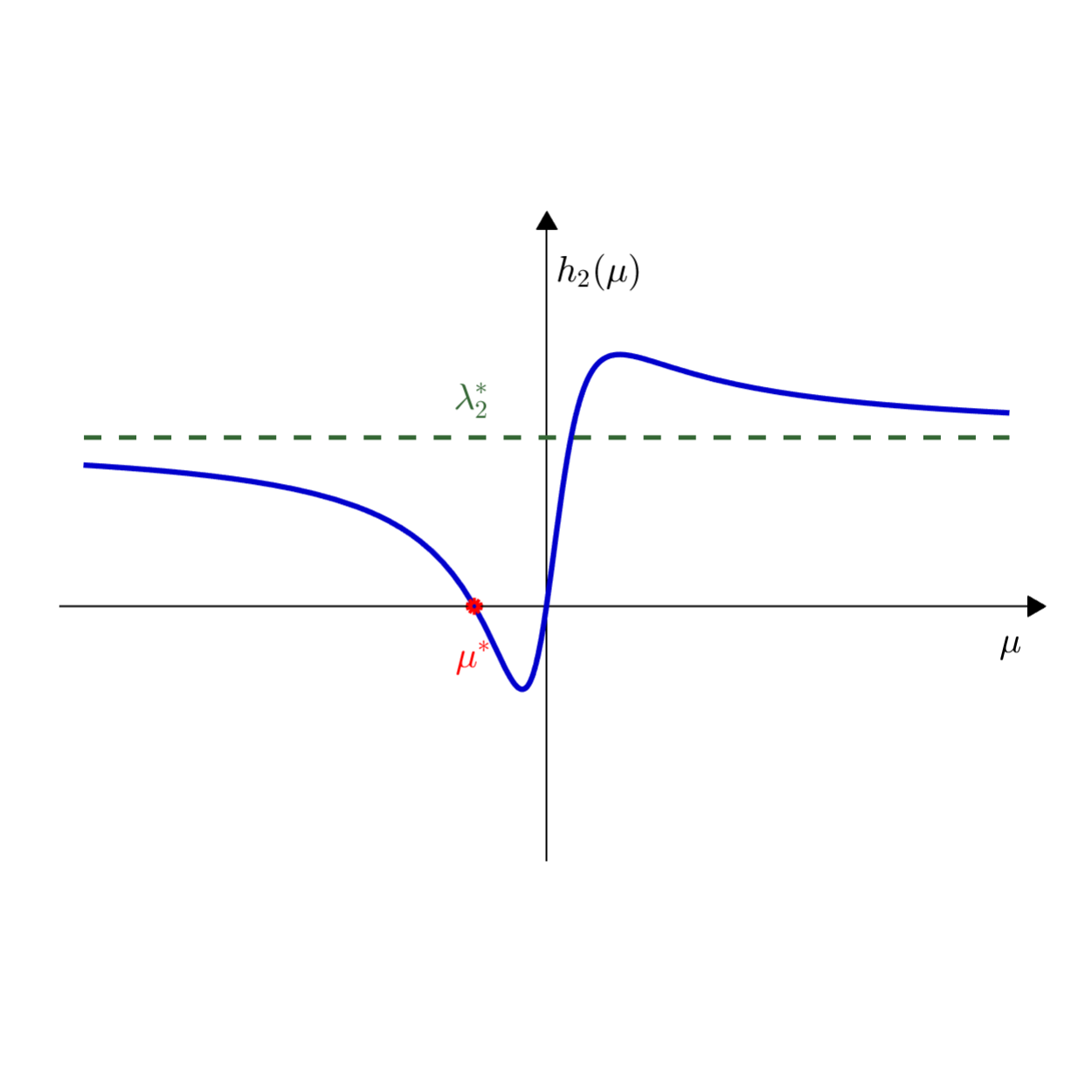}
\caption{Functions $h_1(\m)$ (left) and $h_2(\m)$ (right) in the case $m_1$ and $m_2$ changing sign and $\int_{\O_1}m_1<0$ and $\int_{\O_2}m_2<0$.}
\label{figuram5}
\end{figure}

\begin{figure}[h]
\centering
\includegraphics[scale=.35]{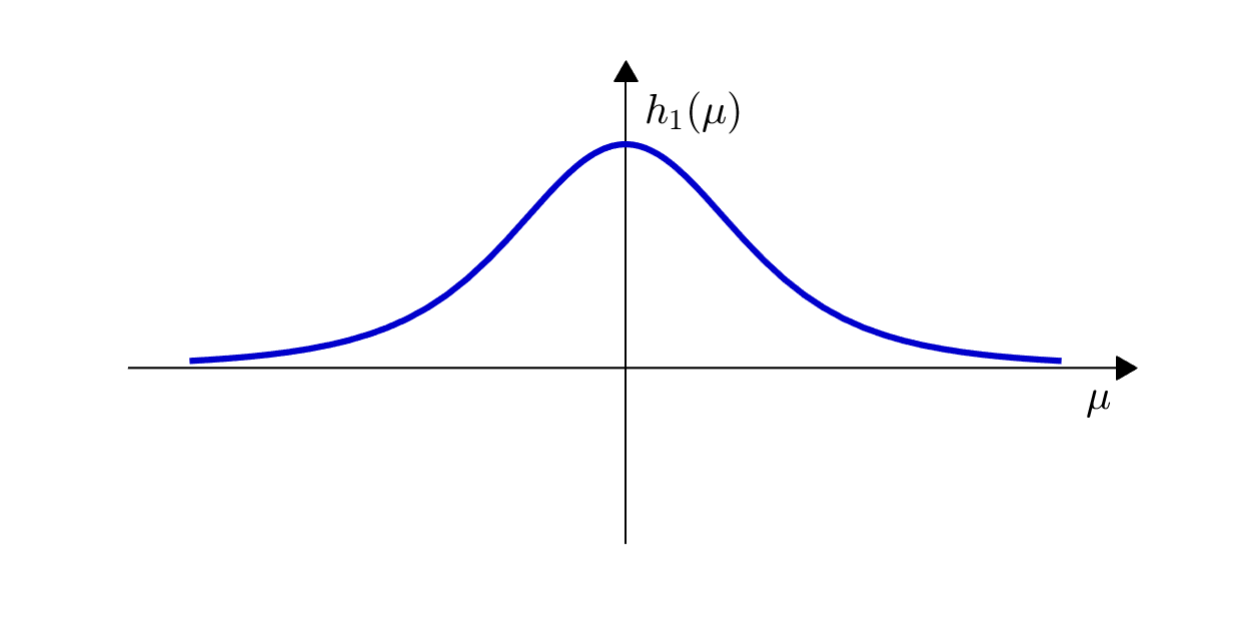}
\includegraphics[scale=.35]{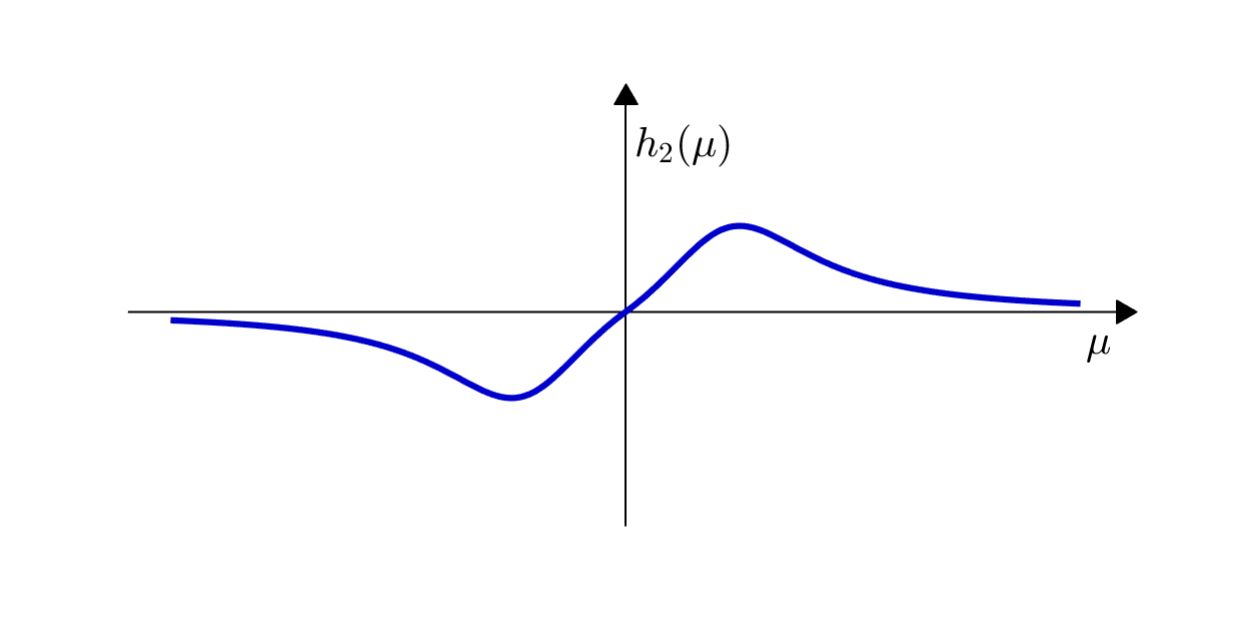}
\caption{Functions $h_1(\m)$ (left) and $h_2(\m)$ (right) in the case $m_1$ and $m_2$ changing sign and $\int_{\O_1}m_1<0$ and $\int_{\O_2}m_2=0$.}
\label{figuram5bis}
\end{figure}

\begin{proposition}
\label{l2cscs}
Assume that $m_i$ changes sign for $i=1,2$ and $
\int_{\O_1}m_1<0,\; \int_{\O_2}m_2<0$. Then, the map $\m\in\R\mapsto h_2(\m)$ is continuous, 
$$
 h_2(\m)
 \left\{
 \begin{array}{ll}
 >0 & \mbox{if $\m<\m^*$ or $\m>0$},\\
 =0 & \mbox{if $\m=\m^*$ and $\m=0$,}\\
 <0 & \mbox{if $\m\in (\m^*,0)$.} 
 \end{array}
  \right.
 $$
Moreover,
$$
\lim_{\m\to \pm\infty}h_2(\m)=\l_2^*.
$$
As a consequence, there exists $\m_{\min}\in (\m^*,0)$ such that
$$
h_2(\m_{\min})=\min_{\m\in\R}h_2(\m)=\l_2^*<0.
$$
\end{proposition}
We have represented in Figures \ref{figuram5} and \ref{figuram5bis} some examples of the maps $h_1(\m)$ and $h_2(\m)$ in the case $m_1$ and $m_2$ changing sign and $\int_{\O_1}m_1<0$ and   $\int_{\O_2}m_2<0$.

\begin{proof}[Proof of Theorem \ref{teoprinci3}:]
\begin{enumerate}
\item By Proposition \ref{l1cscs}, we deduce that
$$
F(\l_1,\l_2)<0\quad\mbox{if $\l_1>\l_1^{\max}$, or $\l_1<\l_1^{\min}$, or $\l_2\geq \l_2^+$ or  $\l_2\leq \l_2^-$}
$$ 
\item Now, we introduce some notation:
$$
h_1(\m):=
\left\{
\begin{array}{ll}
h_1^1(\m) & \mbox{if $\m<\m_{\min}$,}\\
h_1^2(\m) & \mbox{if $\m\in [\m_{\min},\m_{\max}$],}\\
h_1^3(\m) & \mbox{if $\m>\m_{\max}$.}
\end{array}
\right.
$$
%and 
%$$
%h_2(\m):=
%\left\{
%\begin{array}{ll}
%h_2^1(\m) & \mbox{if $\m<\m_*$,}\\
%h_2^2(\m) & \mbox{if $\m>\m_*$.}
%\end{array}
%\right.
%$$

\begin{enumerate}
\item When $\l_1=\l_1^{\max}$, there exists a unique value of $\m$, $\m=\m_{\max}$ such that $h_1(\m_{\max})=\l_1$. The value $h_2(\m_{\max}):=\overline{\l}_2$ verifies that $F(\l_1^{\max},\overline{\l}_2)=0$.

\item Take now $\l_1\in (0,\l_1^{\max})$. Then, there exist $\m^*<\m_2<\m_3$ such that $\l_1=h_1^i(\m_i)$ $i=2,3$, specifically,  
$\l_1=h_1^2(\m_2)=h_1^3(\m_3)$.

Moreover, $\m_2\to \m^*$ and $\m_3\to +\infty$ as $\l_1\to 0$. For these values, $h_2(\m_2)\to 0$ and $h_2(\m_3)\to \l_2^*$. Observe that $h_2(0)=0$. 

\item Consider the case $\l_1\in (\l_1^{\min},0)$. There exists a unique value of $\m_1<\m_2<\m^*$ such that $\l_1=h_1^i(\m_i)$ $i=1,2$, in fact,
$\l_1=h_1^1(\m_1)=h_1^2(\m_2)$.

 In this case, as $\l_1\to 0$, then with $\m_1\to -\infty$ and $\m_2\to \m^*$. Hence,  $h_2(\m_1)\to \l_2^*$ and $h_2(\m_2)\to 0$. 
\item The case $\l_1=\l_1^{\min}$ is analogous to the first case.
\end{enumerate}
Now, we define the maps 
$$
{\cal H}^+(\l_1):=
\left\{
\begin{array}{ll}
h_2((h_1^3)^{-1}(\l_1)) & \mbox{if $\l_1\in (0,\l_1^{\max}],$}\\
h_2((h_1^1)^{-1}(\l_1)) & \mbox{if $\l_1\in [\l_1^{\min},0]$}\\
\end{array}
\right.
$$
and
$$
{\cal H}^-(\l_1):=h_2((h_1^2)^{-1}(\l_1)) \quad\mbox{if $\l_1\in[\l_1^{\min},\l_1^{\max}].$}
$$
\end{enumerate}
This completes the proof.
\end{proof}

\section{Semilinear interface problems}
In this section we study the semilinear problem (\ref{logintro}).
\begin{theorem}
\label{logis}
Problem (\ref{logintro}) possesses a positive solution if and only if $F(\l_1,\l_2)<0$. In case the existence, the positive solution is unique.
\end{theorem}
 \begin{proof}
 Assume that there exists at least a positive solution $(u_1,u_2)$ of (\ref{logintro}). Then, using Proposition~\ref{monoconti} 1.,
 $$
 0=\L_1(-\l_1m_1+u_1^{p_1-1},-\l_2 m_2+u_2^{p_2-1})>\L_1(-\l_1m_1,-\l_2 m_2)=F(\l_1,\l_2),
 $$
% 
% 
% Multiplying by $(\v_1,\v_2)$ a positive eigenfunction associated to $F(\l_1,\l_2)$ we obtain that
% $$
% F(\l_1,\l_2)\left(\int_{\O_1}\v_1u_1+\int_{\O_2}\v_2u_2\right)=-\left(\int_{\O_1} u_1^{p_1}\v_1+\int_{\O_2} u_2^{p_2}\v_2\right),
% $$
 whence we deduce that $  F(\l_1,\l_2)<0$.
 
 On the other hand, assume that $ F(\l_1,\l_2)<0$.  Let $\varphi=(\v_1,\v_2)$ be a positive eigenfunction associated to $F(\l_1,\l_2)$, then
 $$
 \underline{\bu}=(\underline{u}_1,\underline{u}_2)=\varepsilon(\v_1,\v_2),\qquad  \overline{\bu}=(\overline{u}_1,\overline{u}_2)=K(1,1), 
 $$
 it is a pair of sub-supersolution for $\varepsilon$ small and $K$ large. Indeed, $K$ and $\varepsilon$ must verify
 $$
 K^{p_i-1}\geq |\l_i|\|m_i\|_{L^\infty(\O_i)},\qquad
 \varepsilon^{p_i-1} \|\v_i\|_{L^\infty(\O_i)}\leq -F(\l_1,\l_2) \quad i=1,2..
 $$
 Clearly, we can take $ \varepsilon$ small and $K$ large verifying both inequalities and such that $\underline{\bu}\leq \overline{\bu}$ in $\O$.
 
 The uniqueness follows by Theorem 4.3 in \cite{bcs}.
 \end{proof}
 
 \section*{Acknowledgment} MMB, CMR and AS were partially supported by PGC 2018-0983.08-B-I00 (MCI/AEI/FEDER, UE) and by the Consejer\'{i}a de Econom\'{i}a, Conocimiento, Empresas y Universidad de la Junta de Andaluc\'{i}a (US-1380740, P20-01160 and US-1381261). MMB was partially supported by the Consejer\'{i}a de Educaci\'on y Ciencia de la Junta de Andaluc\'{i}a (TIC-0130).

 % \section*{References}

\end{document}